\documentclass[11pt]{article}

\usepackage[usenames]{color}
\usepackage{amssymb,amsmath,latexsym}
\usepackage{amsthm}
\usepackage[mathscr]{eucal}
\usepackage{setspace}
\usepackage{geometry}
 \usepackage{graphicx}
 \usepackage{float}
\usepackage{url}
\usepackage{soul}
\usepackage{cancel}

\newtheorem{lemma}{Lemma}
\newtheorem{theorem}{Theorem}

\newtheorem{proposition}{Proposition}

\newtheorem{remark}{Remark}

\newcommand{\noi}{\noindent}
\newcommand{\beqs}{\begin{equation*}}
\newcommand{\eeqs}{\end{equation*}}
\newcommand{\beq}{\begin{equation}}
\newcommand{\eeq}{\end{equation}}
\newcommand{\beqys}{\begin{eqnarray*}}
\newcommand{\eeqys}{\end{eqnarray*}}
\newcommand{\beqy}{\begin{eqnarray}}
\newcommand{\eeqy}{\end{eqnarray}}

\newcommand{\N}{I\!\!N}

\newcommand{\PP}{I\!\!P}

\newcommand{\E}{I\!\!E}
\newcommand{\R}{I\! \! R}

\newcommand{\pos}{\in(0,\infty)}

\newcommand{\ra}{\rightarrow}
\newcommand{\Ra}{\Rightarrow}

\newcommand{\sig}[2]{\sum_{#1}^{#2}}
\newcommand{\inte}[2]{\int_{#1}^{#2}}

\newcommand{\ab}[3]{{#1}^{#2}_{#3}}

\newcommand{\wt}{\widetilde}
\newcommand{\wh}{\widehat}

\newcommand{\var}{\operatorname{var}}

\mathchardef\mhyphen="2D
\newcommand{\bp}{\begin{pmatrix}}
\newcommand{\ep}{\end{pmatrix}}
\newcommand{\black}[1]{{\textcolor{black}{#1}}{}}
\newcommand{\amy}[1]{{ \textcolor{blue}{(amy says:  #1)}}{}}
\newcommand{\chihoon}[1]{{ \textcolor{blue}{(chihoon says:  #1)}}{}}

\begin{document}

\title{Stationary Distribution Convergence of the Offered Waiting Processes for $GI/GI/1+GI$ Queues in Heavy Traffic}

\author{Chihoon Lee\footnote{E-mail: Chihoon.Lee@stevens.edu}\\ School of Business\\ Stevens Institute of Technology\\ Hoboken, NJ 07030 \\
\and Amy R. Ward\footnote{E-mail: amy.ward@chicagobooth.edu}\\Booth School of Business\\The University of Chicago\\ Chicago, IL 60637 \\
\and Heng-Qing Ye\footnote{E-mail: lgtyehq@polyu.edu.hk}\\Dept. of Logistics and Maritime Studies\\Hong Kong Polytechnic University\\ Hong Kong}
\date{\today}\maketitle

\begin{abstract}
\noi  %We establish the validity of the heavy traffic stationary approximation for a single server queue, operating under the FIFO service discipline, in which each customer abandons the system if his waiting time exceeds his generally-distributed patience time.  \black{This follows from early results of Kingman when the loading factor approaches one from below, but has not been shown in more generality.}  We prove the convergence of \black{the} stationary distributions \black{of the offered waiting time process} and their moments as the traffic intensity approaches unity. %that the hazard rate of the abandonment distribution is scaled and that it is not scaled. %As a consequence, we establish the limit behavior of the stationary abandonment probability {and mean queue-length}.
A result of  Ward and Glynn \cite{WardGlynn2005} asserts that the sequence of scaled offered waiting time processes of the $GI/GI/1+GI$ queue  converges weakly to a reflected Ornstein-Uhlenbeck process (ROU) in the positive real line, as the traffic intensity approaches one. {As a consequence, the stationary distribution of a ROU process, which is a truncated normal, should approximate the scaled stationary distribution of the offered waiting time in a $GI/GI/1+GI$ queue; however, no such result has been proved.} We prove the {aforementioned} convergence, and the convergence of the moments, in heavy traffic, {thus resolving a question left open in~\cite{WardGlynn2005}.  In comparison to Kingman's classical result~\cite{Kingman1961} showing that an exponential distribution approximates the scaled stationary offered waiting time distribution in a $GI/GI/1$ queue in heavy traffic, our result confirms that the addition of customer abandonment has a non-trivial effect on the queue stationary behavior.} %\st{thus  the stationary distribution of ROU provides a valid approximation for the steady-state of the original offered waiting time process. Our study extends Kingman's classical result to incorporate customer abandonments, irrespective of whether the system loading factor approaches 1 from above or below.} %Such a result could follow from early results of Kingman when the traffic intensity approaches one from below, but has not been shown in more generality.
%\footnote{(1) Just write ``Ward and Glynn \cite{WardGlynn2005}'' (remove (2005)?
%(2) The last statement seems ambiguous.
%(Does it mean something like ``our study extends Kingman's classical result to  include customer abandonment''?) -HQ}

\end{abstract}

\noi {\it Keywords: Customer Abandonment; Heavy Traffic; Stationary Distribution Convergence}

\section{Introduction}

There is a long history of studying queueing systems with abandonments, beginning with the early work of Palm~\cite{Palm37} in the late 1930's.  One common objective is to understand the long time asymptotic behavior of such systems, which is governed by the stationary distribution (assuming existence and uniqueness).  However, except in special cases, the models of interest are too complex to analyze directly.  Instead, some researchers have examined the heavy traffic limits of these systems, and developed analytically tractable diffusion approximations (through process-level convergence results).  The question often left open is whether or not the stationary distribution of the diffusion does indeed arise as the heavy traffic limit of the sequence of stationary distributions of the relevant queueing system with abandonment.  Our objective in this paper is to answer this question for one of the most fundamental models, the single server queue, operating under the FIFO service discipline, with generally distributed patience times; that is, the $GI/GI/1+GI$ queue.

Our asymptotic analysis relies heavily on past work that has developed heavy traffic approximations  for the $GI/GI/1+GI$ queue using the offered waiting time process.   The offered waiting time process, introduced in \cite{Baccelli84}, tracks the amount of time an infinitely patient customer must wait for service.  Its heavy traffic limit when the abandonment distribution is left unscaled is a reflected Ornstein-Uhlenbeck process ({see Proposition~\ref{pc}, which re-states a result from \cite{WardGlynn2005} in the setting in this paper}), and its heavy traffic limit when the abandonment distribution is scaled through its hazard rate is a reflected nonlinear diffusion (see \cite{ReedWard2008}).  However, those results are not enough to conclude that the stationary distribution of the offered waiting time process converges, which is the key to establishing the limit behavior of the stationary abandonment probability and mean queue-length.  Those limits were conjectured in~\cite{ReedWard2008}, and shown through simulation to provide good approximations.  However, the proof of those limits was left as an open question. In this paper, we focus on the case when the abandonment distribution is left unscaled (i.e., the heavy traffic scaling studied in \cite{WardGlynn2005}).

{When} the system loading factor is less than one, since the $GI/GI/1+GI$ queue is dominated by the $GI/GI/1$ queue, the much earlier results of \cite{Kingman1961, Kingman1962} for the $GI/GI/1$ queue can be used to establish the weak convergence of the sequence of stationary distributions for the $GI/GI/1+GI$ queue in heavy traffic.  The difficulty arises because, in contrast to the $GI/GI/1$ queue, the $GI/GI/1+GI$ queue can have a stationary distribution when the system loading factor equals or exceeds 1 (see~\cite{Baccelli84}). {\em Our main contribution in this paper is to establish both the convergence of the sequence of stationary distributions and the sequence of stationary moments of the offered waiting time in heavy traffic, irrespective of whether the system loading factor approaches 1 from above or below.}

%{Furthermore, the results in the aforementioned papers only provide the convergence of the stationary distribution, but not its moments.}

An informed reader would recall \cite{Gamarnik:2006} and \cite{BudhLee07}, which establish the validity of the heavy traffic stationary approximation for a generalized Jackson network, without customer abandonment. The proof of the former paper \cite{Gamarnik:2006} relied on certain exponential integrability assumptions on the primitives of the network and as a result a form of exponential ergodicity was established. The latter paper \cite{BudhLee07} provided an alternative proof assuming {the weaker square integrability conditions that are commonly used in heavy traffic analysis}.  Our analysis is inspired by the methodology developed in the latter work \cite{BudhLee07}. However, the main difficulty in extending their methodology to the current model is that the known regulator mapping under customer abandonment is only \emph{locally} Lipschitz (that is, the Lipschitz constant depends on time parameter) whereas the proofs of both \cite{Gamarnik:2006} and \cite{BudhLee07} critically rely upon the global Lipschitz property of the associated regulator mapping. {More precisely, the approaches in \cite{Gamarnik:2006} and \cite{BudhLee07} make use of
the global Lipschitz continuity property of an associated regulator (Skorokhod) mapping to help
convert the given moment bound of primitives (the inter-arrival and
service times) to the bound of the key performance measures
(the waiting time or queue-length processes).  Such a property is
not available for the model under study.}

%we require a different regulator mapping than the one used in their paper, that depends on the customer abandonment distribution;
In connection to the aforementioned technical issue, the studies of \cite{YeYao12-interchange, YeYao18} extend
the works of \cite{Gamarnik:2006} and \cite{BudhLee07} to a wider range
of stochastic processing networks, e.g., the multiclass queueing network
and the resource-sharing network, by relaxing the requirement of the {afore}mentioned
Lipschitz continuity.
However, their study in \cite{YeYao12-interchange, YeYao18} deals with networks that have
heavy-traffic limits satisfying the {\it linear} dynamic complementarity problem,
i.e., the state process depends on the ``free process'' (and the regulating
process as well) linearly.
Therefore, their results do not apply to our $GI/GI/1+GI$ model directly,
as the resulting heavy-traffic limit is a reflected Ornstein-Uhlenbeck process
and the state process of this limit (i.e., $V(\cdot)$ in \eqref{3.58} below)
depends on the ``free process''  (i.e., the drifted Brownian motion in \eqref{3.58}) in a  nonlinear manner.
Nevertheless, their hydrodynamics approach is adapted to establish a key property,
i.e., the uniform moment stability of the offered waiting time process
(see Section \ref{section:UniformMoment}), in our paper.
%\red{[1. Can we say more about why the reflected O-U process does not satisfy the linear dynamic complementarity problem.  I am guessing this is because there is a state space dependence in the drift that is not present in the other work?]} \red{[2. How [11] and [4] get around the need for global Lipshitz continuity?  What do you think about adding a sentence on the high level approach in [11] and [4]?]}

%In order to overcome such a technical restriction, we adapt the relevant hydrodynamical scaling analysis from the works of \blue{(Brief remark on the current hydrodynamical scaling approach.)}
{A closely related paper is that of  Huang and Gurvich  \cite{HuangGurvich18}, which studies the Poisson arrival case (i.e., $M/GI/1+GI$ queue) and shows the associated Brownian model is accurate uniformly over a family of patience distributions and universally in the heavy-traffic regime. For instance, Section EC.3.1 therein corresponds to the critically loaded regime, as considered in this paper. Their approach is based on the generator comparison methodology, and owing to the Poisson arrivals, it is enough to consider a one-dimensional process with a simple generator, whereas with general arrival processes, one needs to consider a two-dimensional process (tracking, e.g., the residual arrival times) and correspondingly more  complicated generator.}

In comparison to results for many-server queues, the process-level convergence result for the $GI/GI/N+GI$ queue in the quality-and-efficiency-driven regime was established in \cite{MandelbaumMomcilovic2012} when the hazard rate is not scaled, and in \cite{ReedTezcan2012}, under the assumption of exponential service times,  when the hazard rate is scaled.  Neither paper establishes the convergence of the stationary distributions.  That convergence is shown under the assumption that the abandonment distribution is exponential and the service time distribution is phase type in \cite{DaDiGa2014}.  The question remains open for the fully general  $GI/GI/N+GI$  setting.  There has been some progress made in this direction in \cite{KR2012}, which establishes the convergence of the sequence of stationary distributions under fluid scaling in the aforementioned fully general $GI/GI/N+GI$ setting.

The remainder of this paper is organized as follows. We conclude this section with a summary of our mathematical notation.  In Section \ref{subsection:modelassumption}, we set up the model assumptions and recall the known process-level convergence results for the $GI/GI/1+GI$ queue.  In Section \ref{section:SSmoments}, we state our main result, that gives the convergence of the stationary distribution of the offered waiting time process, and its moments.    %Section \ref{section:UniformMoment} shows how to obtain bounds on the moments of the \black{scaled state process} that are uniform in the heavy traffic scaling parameter ($n$).
{To prove the main result, we first obtain bounds on the moments of the scaled state process that are uniform in the heavy traffic scaling parameter ($n$) in Section \ref{section:UniformMoment}. The proofs of lemmas in this section are technically involved, and are delayed to Section \ref{section:LemmaProof}.  Lastly, we use uniform moment bounds established in Section \ref{section:UniformMoment} to prove
our main result in Section \ref{section:TheoremProof}.}

%Section \ref{section:LemmaProof} presents the proof of all lemmas. Lastly, we use the established uniform moment bounds in Section \ref{section:UniformMoment}  to prove our main result in Section \ref{section:TheoremProof}.

{\bf Notation and Terminology.} Use the symbol ``$\equiv$'' to stand for equality by definition. The set of positive integers is denoted by $\N$ and denote $\N_0\equiv \N\cup\{0\}$.  Let $\R$ represent the real numbers $(-\infty, \infty)$ and $\R_+$ the non-negative real line $[0,\infty)$.  For $x,y \in \R$, $x \vee y \equiv \max\{x,y\}$ and $x \wedge y \equiv \min\{x,y\}$.   The function $e(\cdot)$ represents the identity map; that is, $e(t) = t$ for all $t \in \R_+$. For  $t\in \R_+$ and a real-valued function $f$, define $\| f \|_t \equiv \sup_{0 \leq s \leq t} |f(s)|$. Let $D(\R)\equiv D(\R_+, \R) $ be the space of right-continuous functions $f:\R_+\ra\R$ with left limits, endowed with the Skorokhod $J_1$-topology (see, for example,~\cite{Billingsley:1999}). Lastly, the symbol ``$\Ra$'' stands for the weak convergence; we make this explicit for stochastic processes in $D(\R)$, otherwise, it is used for weak convergence for a sequence of random variables. % vs weak convergence for processes in $D(\R)$. %Define $D(\R_+,X)$ and ``$\Ra$''. %Similarly, define $D(\R_+) \equiv D(\R_+, \R_+)$ and $D^2(\R_+)\equiv D(\R_+, \R_+)\times D(\mathbb R_+, \R_+)$  \blue{Chihoon, you mean $D^2(\R_+)\equiv D(\R_+, \R_+)\times D(\R_+, \R_+)$, right?  At least, that is what we need when we define GRM. }

%\blue{Amy:  This subsection still needs work.  Need to double-check that everything aligns with the main body.  Also, in the second paragraph we need to differentiate between weak convergence for a sequence of random variables vs weak convergence for processes in $D(\R)$.  This is important because our main result is not a process-level convergence result.}

\section{The Model and \black{Known} Results} \label{subsection:modelassumption}
The $GI/GI/1\textcolor{black}{+GI}$ model having FIFO service is built from three independent i.i.d. sequences of nonnegative random variables {$\{ u_i, i \geq 2\}$, $\{v_i,i\geq 1\}$, $\{d_i,i\geq 1\}$,  that are representing inter-arrival times, service times, and patience times, respectively, and are defined on a common probability space $(\Omega, \mathcal{F}, \PP)$.}   At time 0, the previous arrival to the system occurred at time $t_0^n <0$, so that $|t_0^n|$
represents the time elapsed since the last arrival in the $n$-th system.  We let $u_1$ be the random variable representing the remaining time conditioned on $|t_0^n|$ time units having passed; that is,
\[
\PP(u_1 >x) = \PP(u_2>x | u_2 > |t_0^n| ).
\]
We let $F(\cdot)$ represent the distribution function associated with the patience time $d_1$ and, {consistent with \cite{Baccelli84}, assume $F(\cdot)$ is proper.}
The system primitives are assumed to satisfy:
\begin{enumerate}
\item[($\mathbb A$1)] {For some $p\in (2,\infty)$, {$\E[u_2^{p}+v_2^{p}]<\infty$,}} and $F'(0)\in(0,\infty)$.
%\blue{[Chihoon: We need $2+\delta$ moment condition for the stationary distribution convergence, and slightly higher moment condition for the stationary moment convergence.The latter is needed for the U.I. to hold. Let us think how we can revise/streamline the associated assumptions . (The $\delta$ was used in the proof of Lemma 3. )]}
\end{enumerate}

%\hq{We may fix $p\ge2$ here, and don't touch it anymore
%below. Then, we don't need mention ``p=2'' in Thm 1 and 2.
%(Average readers won't care about the differences.
%And we may comment like ``it suffices to assume ($\mathbb A$1) with $p=2$
%for Thm 1 and Thm 2(a)''
%following the statement of Thm 2).
%We don't need to mention $p\ge 1$ in Prop 1,
%and simply say like ``Prop 1 implies (56) via Jensen inequality''
%at the beginning of Sec 6.
%(I think we don't even have to mention ``Jensen''.)
%}

We consider a sequence of systems indexed by $n\geq1$ in which the arrival rates become large and service times small. By  convention, we use superscript $n$ for any processes or quantities associated with the $n$-th system.  The arrival and service rates in the $n$-th system are $\lambda^n$ and $\mu^n$ and satisfy the following heavy traffic assumption:
\begin{enumerate}
\item[($\mathbb A$2)] $\lambda^n\equiv n\lambda$, $\lim_{n \rightarrow \infty} \frac{\mu^n}{n} = \lambda \in (0,\infty)$ and $\lim_{n \rightarrow \infty} \sqrt{n} \left(\lambda - \frac{\mu^n}{n} \right) = \theta \in \R.$
\end{enumerate}
The $i$-th arrival to the $n$-th system occurs at time
\[
t_i^n \equiv \sum_{j=1}^i \frac{u_j}{\lambda^n} \,\,\mbox{ for }\,\, i \in\N,
\]
{and} 
has service time
\[
    v_i^n \equiv \frac{v_i}{\mu^n}\,\, \mbox{ for }\,\, i \in \N,
\]
and abandons without receiving service if processing does not begin by time $t_i^n + d_i$.

\noi \textbf{The Offered Waiting Time Process}

The offered waiting time process, first given in~\cite{Baccelli84}, tracks the amount of time an incoming customer at time $t$ has to wait for service.  That time depends only upon the service times of the non-abandoning customers already waiting in the queue, that is, those waiting customers whose patience time upon arrival exceeds their waiting time.  For {$t\geq0$}, the \emph{offered waiting time} process having initial state $V^n(0)$  has the evolution equation \beq \label{2.1.1}  V^n(t)=  \black{V^n(0) } +  \sig{j=1}{A^n(t)} v^n_j\mathbf{1}_{[V^n(t^n_j-)<d_j]} - \inte{0}{t}\mathbf{1}_{[V^n(s)>0]}ds\geq0, \eeq
 where \beq\label{254} A^n(t) \equiv \max\{ i \in \N_0 : t_i^n \leq t \} \eeq
 is a delayed renewal process when $|t_0^n| >0$ and is a regular (non-delayed) renewal process when $t_0^n = 0$.
 The quantity $V^n(t)$ can also be interpreted as the time needed to empty the system from time $t$ onwards if there are no arrivals after time $t$, and hence it is also known as the workload at time $t$.  {The initial state $V^n(0)$ is 0 if no job is in service and otherwise represents the total workload of all jobs that arrived prior to time 0 and that will not abandon before their service begins.}

\noi {\bf Reflected Ornstein-Uhlenbeck Approximation.}

We consider the one-dimensional reflected Ornstein-Uhlenbeck process $V\equiv \{V(t)\}_{t\geq 0}$
%Let $(V,L)$ be the unique solution to the reflected stochastic differential equation
 \beq\label{3.58}
 \begin{array}{l}
V(t)=V(0) +   \sigma W(t)+ {\frac{\theta}{\lambda}} t- F'(0)\inte{0}{t}V(s)ds + L(t) \geq 0 \\
 \mbox{subject to:  } L \mbox{ is non-decreasing, has } L(0)=0 \mbox{ and } \int_{0}^{\infty}V(s)dL(s)=0,
 \end{array} \eeq
where $\{W(t):t\geq0\}$ denotes a one-dimensional standard Brownian motion, and the infinitesimal variance parameter is
 \[
 {\sigma^2  \equiv \lambda^{-1}( \var(u_2) + \var(v_2))}.
  \]
{Given $\mathcal F_0$-measurable initial condition $V (0)\geq 0$, the strong existence and pathwise uniqueness of solution $(V,L)$ to the  stochastic differential equation \eqref{3.58} hold for the data $(V (0), W)$, i.e., the solution is adapted to $(\mathcal F^W_t \vee \mathcal F_0)_{t\geq0}$ (see, e.g., \cite{Zhang94}).}

%When the initial condition satisfies $\sqrt{n} V^n(0) \Rightarrow V(0)$ as $n\ra\infty$, Theorem 1(a) of  \cite{WardGlynn2005} adapted to the setting in this paper  (cf. Remark \ref{ds1} below) shows that
{The following weak convergence result is a simple modification of Theorem 1(a) of  Ward and Glynn \cite{WardGlynn2005} and we provide its proof in the Appendix for the sake of completeness.
\begin{proposition}\label{pc}
Assuming $\sqrt{n} V^n(0) \Rightarrow V(0)$ as $n\ra\infty$, we have
  \begin{equation} \label{eq:Vapprox}
 \sqrt{n}V^n \Ra V \,\,\mbox{ in }\,\, D(\R) \,\,\mbox{ as } \,\,  n \ra \infty.
 \end{equation}
 \end{proposition}}
%where $D(\R)$ is the space of right-continuous functions $f:\R_+\ra\R$ with left limits, endowed with the Skorokhod $J_1$-topology (see, for example,~\cite{Billingsley:1999}), and the symbol ``$\Ra$'' stands for the weak convergence.

%\blue{We should not show the convergence result for $Q$ above, or mention it below, if we cannot guarantee the existence of a stationary distribution for the queue-length.}

\section{{The Stationary Distribution Existence and Convergence Results}} \label{section:SSmoments}
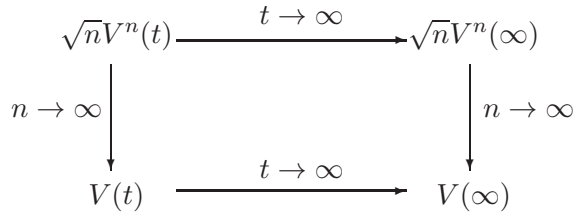
\begin{figure}[!t]
\centering
\unitlength 1mm % = 2.85pt
\linethickness{0.4pt}
\ifx\plotpoint\undefined\newsavebox{\plotpoint}\fi % GNUPLOT compatibility
%\begin{picture}(114.75,145.25)(0,0)
%\begin{picture}(75.75,57)(50,100)
\begin{picture}(75.75,31)(50,115)
\put(64,139.5){\makebox(0,0)[cc]{$\sqrt{n}V^n(t)$}}
\put(111.5,139.5){\makebox(0,0)[cc]{$\sqrt n V^n(\infty)$}}
\put(64,118){\makebox(0,0)[cc]{$V(t)$}}
\put(111.5,118){\makebox(0,0)[cc]{$V(\infty)$}}
\put(72,138.75){\vector(1,0){30.5}}
\put(72,118.75){\vector(1,0){30.5}}
\put(63.25,135.5){\vector(0,-1){14}}
\put(111,135.5){\vector(0,-1){14}} \put(83,141){$t\rightarrow
\infty$} \put(83,120.5){$t\rightarrow \infty$}
\put(50,128.5){${n\rightarrow\infty}$}
\put(112.75,128.5){${n\rightarrow\infty}$}
\end{picture}
\caption{\em A graphical representation of the limit interchange. } \label{120}
\end{figure}

{The weak convergence (\ref{eq:Vapprox}) motivates approximating the scaled stationary distributions for $V^n$, and its moments, using the stationary distribution of $V$, and its moments.  This requires establishing the limit interchange depicted graphically in Figure~\ref{120}.  When the limit $n\rightarrow \infty$ is taken first, and the limit $t \rightarrow \infty$ is taken second, the convergence is known.  More specifically, the convergence as $n \rightarrow \infty$ was established in (\ref{eq:Vapprox}), and Proposition 1 in~\cite{WardGlynn-ROU2003} shows
\[
V(t) \Rightarrow V(\infty) \,\mbox{ as }\, t \rightarrow \infty
\]
for $V(\infty)$ a random variable having density
\begin{equation} \label{eq:VstationaryDensity}
 f(x) = \frac{b^{-1}\phi(\frac{x-m}{b})}{1-\Phi(-\frac{m}{b})}\, \mbox{ for }\, x \geq 0,
\end{equation} where $m\equiv \theta/(\lambda F'(0))$, $b\equiv \sigma/{\sqrt{2F'(0)}}$ and $\phi(\cdot)$, $\Phi(\cdot)$ denote the pdf and cdf of the standard normal distribution, respectively. In other words, $V(\infty)$ is distributed as truncated normal with mean $m$ and variance $b^{2}$, conditioned to be on $\R_+$.}

In this Section, we state our main results, first that a unique stationary distribution exists for each system $n$, and second that
the convergence in Figure~\ref{120} is valid when the limit is first taken as $t \rightarrow \infty$ and second taken as $n \rightarrow \infty$.  In order to do this, we  first specify the relevant Markov process.

The offered waiting time process $\{V^n(t):t\geq0\}$ alone is not Markovian due to the {remaining} arrival time. {(In contrast, the offered waiting time process tracked only at customer arrival times is a Markov chain with state space $\R_+$; see~\cite{Baccelli84} and the recursive equations therein.)}
% However, we can augment the state space to define a vector-valued process that is not only Markov, but is also strong Markov.  To do this, we follow~\cite{Davis84}, and, for convenience to the reader, we use notation consistent with that paper.
Defining the remaining arrival time (i.e., the forward recurrence time of the arrival process)
\[
\tau^n (t)\equiv t_{j+1}^n - t \,\mbox{ for }\, t \in [t_j^{{n}}, t_{j+1}^{{n}}), \,\,j \in \N_0,
\]
where {$\tau^n(0) =u_1$},  the vector-valued process $$\mathbb X^n \equiv \{(\tau^n(t), V^n(t)):t\geq0\}$$ having state space $\mathbb S \equiv  \R_+ \times \R_+$ is strong Markov (cf. Problem 3.2 of Chapter X in \cite{Asmussen2003} and also Section 2 in  \cite{Davis84}). %The queue has two states, ``busy''  and ``empty'', and will correspond to the states $\nu^n(t)= 1,0$, respectively. % indicate one of the two states of the system, ``busy''  and ``empty''. %, and these will   Let $\nu^n(t) \in \{0,1\}$ indicate whether or not the queue is empty, and assume $\nu^n(0) = 0$ (representing an empty queue) \red{if and only if $V^n(0) = 0$.}

%\begin{lemma} \label{lemma:strongMarkov}
%For each $n\geq1$, $\mathbb X^n \equiv \{(\nu^n(t), \tau^n(t), V^n(t)):t\geq0\}$ having state space $\mathbb S \equiv \{0,1\} \times \R_+ \times \R_+$ is a strong Markov process.
%\end{lemma}
%\noindent The proof of Lemma~\ref{lemma:strongMarkov} is presented at the end of this Section.

Ensuring the existence of a stationary distribution requires the following {technical condition} on the interarrival time.
More precisely, this assumption is used {in the proof of Theorem \ref{stability} below to verify a petite set requirement that implies positive Harris recurrence (see Lemma~\ref{petite-set}), from which the existence of  a unique stationary distribution follows immediately.}
{Such an assumption}  has been frequently used in the literature; for example, Proposition 4.8 in \cite{Bramson:2008}, Lemma 3.7 in \cite{MeynDown1994}, and Theorem 3.1 in \cite{Dai:1995}. %, see also Chapter VII of \cite{Asmussen2003}. % (see the proof of Theorem \ref{stability}).
For $x=(\tau, v) \in \mathbb S$, define its norm $|x|$ as $|x|\equiv\tau+v$. % To state that technical condition, we
Define the norm $|\mathbb X^n(t)|$ of $\mathbb X^n(t)$ to be the sum of the offered waiting time and the remaining arrival time at $t$, that is, \[|\mathbb X^n(t)| \equiv V^n(t) + \tau^n(t), \,\, t\geq0\,. \]
\begin{enumerate}
\item[($\mathbb A$3)] The i.i.d. interarrival times $\{u_i, i\geq2\}$ are unbounded (that is, $\PP(u_2\geq u)>0$ for any $u>0$). %and have a positive continuous density function}.} %there exist some integer $k_0\geq2$ and some function $p(x)>0$ on $\R_+$ with $\inte{0}{\infty} p(x) dx>0$ such that $\PP(u_2\geq x)>0$ for any $x>0$ and $\PP(a\leq u_2+\cdots+u_{k_0}\leq b) \geq \inte{a}{b} p(x)dx$ for any $0\leq a<b$.} % $\nu\{[L,\infty)\}>0$ for all $L>0$. Also, $\nu$ is spread out, that is, for some $k_0\geq2$, the $(k_0-1)$-fold convolution $\nu^{(k_0-1)*}$ is nonsingular with respect to Lebesgue measure.} %\textcolor{red}{Delete this: There exists a probability distribution $a$ on $[0,\infty)$ such that $$\inte{0}{\infty} \bar F_A(t)(1+F_S(t)\bar F_D(t))a(dt)<\infty,$$ where $\bar F_A$ denotes the complementary cdf (ccdf) of interarrival times $u_2,u_3, \ldots$, $F_S$ the cdf of service times $v_1, v_2, \ldots$, and $\bar F_D$ the ccdf of patience times $d_1,d_2, \ldots$.} %\textcolor{red}{The interarrival time distribution $\nu$ is unbounded, so that $\nu\{[L,\infty)\}>0$ for all $L>0$. Also, $\nu$ is spread out, that is, for some $k_0\geq2$, the $(k_0-1)$-fold convolution $\nu^{(k_0-1)*}$ is nonsingular with respect to Lebesgue measure.}
\end{enumerate}
%\noindent   We refer the reader to Lemma 3.2 of \cite{Dai:1995} for a related sufficient conditions for ($\mathbb A$3) (i.e., the interarrival times are unbounded and spreadout) and also \cite{Meyn:Tweedie:1993[3]} for the concept and role of the petite set for stability analysis of continuous time Markov processes on a general state space.  \blue{(Modify this later.)}  %. The interarrival times are unbounded and spreadout,

%\begin{remark}
%\textcolor{red}{%It is evident that the assumption $(\mathbb A3)$ holds when the interarrival, service, patience times are all exponential random variables.
%This assumption is used when verifying a petite set requirement for the sought-after positive Harris recurrence property in the following theorem.}
%\end{remark}

\begin{theorem}\label{stability} \textsc{(Stationary Distribution Existence)} {Assume $(\mathbb A1)$--$(\mathbb A3)$.}
For each $n\geq1$,
there exists a unique stationary probability distribution for the Markov process ${\mathbb X}^n$.
\end{theorem}

Now that we know the stationary distribution, denoted as $\pi^n$, exists for each fixed $n$, we can establish its convergence to the stationary distribution of the diffusion (\ref{3.58}) given in (\ref{eq:VstationaryDensity}).  To state this result, we require the diffusion-scaled process:
\[  \widetilde{\mathbb X}^n \equiv (\wt\tau^n, \wt{V}^n), \mbox{ where } \wt{\tau}^n(t)\equiv \sqrt n \tau^n(t) ,\,\,\, \wt{V}^n(t)\equiv \sqrt n V^n(t)\,. \]
Notice %the component $V^n$ scales with $n$ but the other components do not.
the time is not scaled in the process $\wt{V}^n$ because the arrival and service rate parameters are scaled instead (from $(\mathbb A2)$, both $\lambda^n$ and $\mu^n$ are order $n$ quantities); therefore, scaling the state by $\sqrt{n}$ produces the traditional diffusion scaling. {Also, a motivation behind the scaling for the residual arrival time $\tau^n$ comes from the way how the $\sqrt n$ diffusion scaling affects the $\tau^n$  under the arrival rate $\lambda^n=n\lambda$ (cf. see \eqref{ds} below)}.
%the fact that we have obtained the diffusion scaled process $\tilde{A}^n$ by scaling the rate $\lambda^n$, which then affects the initial time elapsed since the last arrival.)} %Since time is not scaled, there is no need to scale the {remaining} arrival process $\tau^n$.
%The component $\nu^n$ does not need to scale because it is an indicator function for the queue being busy or empty.

%\green{In the below, should $\mathbb S$ by $\widetilde{\mathbb X}^n$? }
\begin{theorem}\label{convstat} \textsc{(Stationary Convergence)}
%Assume $(\mathbb A1)$--$(\mathbb A3)$.  For any fixed large enough $n$,
Assume $(\mathbb A1)$--$(\mathbb A3)$.
\begin{enumerate}
\item[(a)] \textsc{(Distribution)}    %Under the same assumption of Theorem \ref{stability},
Denote by $\pi^n_0$ the marginal distribution of $\pi^n$ on the {second} coordinate of {$\widetilde{\mathbb X}^n$,}  i.e., $\pi^n_0(A)= \pi^n(\R_+ \times A)$ for  $A\in\mathcal B(\R_+)$. Let $\wt{V}^n(\infty)$ be a random variable having distribution $\pi^n_0$ and also $V(\infty)$ a random variable having density \eqref{eq:VstationaryDensity}.
We have that $\widetilde{V}^n(\infty) \Rightarrow V(\infty)$ as $n \rightarrow \infty$.
\item[(b)] \textsc{(Moments)}
%If also $\E [u_1+ v_1]^{q} < \infty\, \mbox{ for some }\, q\in (2,\infty)$,
%Assume $(\mathbb A1)$--$(\mathbb A3)$. %and let $q\equiv p+\delta \in(2,\infty)$.
%  For any  $m\in(0,p-\delta-1)$ with $\delta\in(0,p-2)$,
{For any  $m\in(0,p-1)$,} % with $q\in[2,p)$
\begin{equation*} \E[(\wt V^n(\infty))^m] \ra \E[(V(\infty))^m] \,\mbox{ as }\, n \rightarrow \infty.  \end{equation*}
\end{enumerate}
\end{theorem}

\begin{remark} \textsc{(Queue-length Convergence)}  The queue-length process $Q^n(t)$ represents the number of customers that are in system at time $t >0$, either waiting or with the server.  In contrast to $V^n$, $Q^n$ includes customers that will eventually abandon but have not yet done so.  When the initial condition satisfies $Q^n(0) / \sqrt{n} \Rightarrow Q(0)$ as $n \rightarrow \infty$, Theorem 3 in \cite{WardGlynn2005} shows that
\[
    \frac{Q^n}{\sqrt{n}} \Rightarrow \lambda V \,\,\mbox{ in }\,\, D(\R) \,\,\mbox{ as } \,\,  n \ra \infty.
\]
In other words, recalling \eqref{eq:Vapprox}, a process-level version of Little's law holds.  This suggests that a version of Theorem~\ref{convstat} should hold for the queue-length process as well.  However, proving this is more involved technically due to the need to track {\em both} customers in queue that will eventually receive service {\em and} customer in queue that will eventually abandon; see the ``potential queue measure'', a measure-valued state descriptor in Section 2.2 of  \cite{KR2012}, {and also see Figure 4 in~\cite{PuhaWard2019} for a graphic depiction of that measure}.  This is the reason we leave that analysis as future research.
\end{remark}

%Finally, to end this section, we provide the details to establish that ${\mathbb X}^n$ is a strong Markov process, as asserted earlier.

\section{{Uniform Moment Estimates}} \label{section:UniformMoment}

The proofs of both Theorems~\ref{stability} and~\ref{convstat} rely on a tightness result for the family of stationary distributions  of {$\{\wt{\mathbb X}^n\}_{n\geq1}$}. The key to the desired tightness is to obtain \emph{uniform} (in $n$) bounds for the  moments  of the stationary distributions. Henceforth, we use the subscript $x$ to denote the scaled Markov process $\wt{\mathbb X}^n$ has an initial state  $(\wt\tau^n(0), \wt V^n(0))= (\tau, v)\equiv x\in\mathbb S$.  Our convention is to subscript any process that depends on the initial state $x \in \mathbb S$ by $x$.  Then, $\wt{V}_x^n$ has initial state $V^n(0) = v/\sqrt{n}$ and is defined from the process $V_x^n$ in (\ref{2.1.1}) that uses the (delayed) renewal process $A_x^n$ in (\ref{254}). % having $|t_0^n| = \tau$ for all $n$.
Recall the norm $|x|$ of $x=(\tau,v)\in \mathbb S$ is defined as
$ |x| \equiv \tau + v. $% \,\,\mbox{ for }\,\, x \in \mathbb S.
%\blue{(Chihoon, I would like to eliminate the next sentence because it is repetitive.)}Proposition \ref{moment1} below provides estimates on moments of the (scaled) state process that are uniform in the scaling parameter $n$.

%In what follows, if the initial condition of the Markov process $\mathbb X^n$ is $x$ for some $x\in\mathbb S$, i.e., $(\nu^n(0), \tau^n(0), V^n(0))= (\nu, \tau, v)\equiv x\in\mathbb S$,

%\red{(Chihoon, do we need $(\mathbb A3)$ to prove the below Proposition?  I think we only need $(\mathbb A1)$ and $(\mathbb A2)$. )}
\begin{proposition}\label{moment1}
Assume $(\mathbb A1)$--$(\mathbb A2)$.  {Let $q\in[1,p).$}
%\blue{[We would be using different symbol than $p$ here.]}
There exists  $t_0\pos$ such that for all $t\geq t_0$,
\beq\label{2.1} \lim_{|x|\ra\infty}\sup_{n}\frac{1}{|x|^{{q}}}\E\left[|\wt{\mathbb X}^n_x(t|x|)|^{{q}} \right]=0. %, \quad p\in(0,q],
\eeq
%{where $q$ is as in Theorem \ref{convstat}.}
\end{proposition}

{Before proving Proposition \ref{moment1}, we provide a  roadmap about how it will be used in the proofs of the main results.  Proposition \ref{moment1} yields uniform (in the scaling parameter) moment bounds for the scaled Markov  process of the $GI/GI/1+GI$ system. We use such uniform moment estimates, in conjunction with the Lyapunov function methods of Meyn and Tweedie \cite{Meyn:Tweedie:1993[3]} and Dai and Meyn \cite{Dai:Meyn:1995}, in order to obtain time uniform moment bounds (via weighted return time estimates) for  {the aforementioned scaled Markov process}. Then, the sought-after moment bounds for the stationary distributions, uniform in the scaling parameter, readily follow (see \eqref{eq:to-show} and \eqref{eq:needed1} in the proof of Theorem \ref{convstat}) and this yields tightness of the collection of the stationary distributions and hence Theorem \ref{convstat}.}  %Indeed, the proof of Theorem \ref{convstat} is based on such key estimates on moments of the state process that are uniform in both time and the scaling parameter . }

%\red{(Slightly modify this later.)} {The proof of Proposition~\ref{moment1} relies on a martingale representation of the offered waiting time process and the continuity properties of a nonlinear generalized regulator mapping.  We first provide this setup, and second give the proof.} In doing so, we consider  the case that the hazard rate is scaled (that is, $(\mathbb A1)$--$(\mathbb A4)$ hold). The case that the hazard rate is not scaled is covered by minor adjustments to the proofs, which are shown at the end of proof of Theorem~\ref{convstat}.

{The crux in proving Proposition~\ref{moment1} lies in two versions of pathwise stability results
(Lemmas \ref{lem-10a} and \ref{lem-10b}), whose intuitive ideas are provided right after stating those results.} The proof of Proposition~\ref{moment1} relies on a martingale representation of the offered waiting time process. % and a Lipschitz continuity property of an $n$-dependent regulator mapping.
We first provide this setup, and second give the proof.  %\blue{(Chihoon, I eliminated much of the previous paragraph, because the reader will have a hard time understanding.  For example, $\wt{\epsilon}_x^n$ was referred to in that paragraph but had not been defined.)}

\subsection{Martingale Representation and Diffusion Scaling} %of the Offered Waiting Time Process}
\label{subsubsection:MGrepresentation}

Define the $\sigma$-fields $(\wh{\mathcal F}^n_i)_{i\geq1}$ where \[  \wh{\mathcal F}^n_i\equiv \sigma((t^n_1,v^n_1,d_1), \ldots, (t^n_i,v^n_i,d_i), t^n_{i+1})\subseteq \mathcal F, \] and let $\wh{\mathcal F}_0^n\equiv\sigma(t_1^n)$. %Next, we define the associated continuous time filtration $({\mathcal F}_t^n)_{t\geq0}$ by ${\mathcal F}_t^n\equiv \wh{\mathcal F}_{[nt]}^n$.
Notice that $ V^n_x(t^n_i-)$ is $\wh{\mathcal F}^n_{i-1}$-measurable and the patience time  $d_i$ of the $i$-th customer is independent of $\wh{\mathcal F}^n_{i-1}.$  Hence, \beq\label{2.3.1} \PP[V^n_x(t^n_i-)\geq d_i|\wh{\mathcal F}^n_{i-1}]=F(V^n_x(t^n_i-)), \quad i=1,2, \ldots,  \eeq holds almost surely, recalling that  $F$ is the distribution function of $d_i$. We then have a martingale with respect to the filtration $(\wh{\mathcal F}^n_i)_{i\geq1}$ given by
\[  M^n_{x}(i)\equiv  \sig{j=1}{i} \left(\mathbf{1}_{[V^n_x(t_j^n-)\geq d_j]}  - \E (\mathbf{1}_{[V^n_x(t_j^n-)\geq d_j]} |\wh{\mathcal F}_{j-1}^n)\right).
\]
Using \eqref{2.3.1}, we also see that for all $i\in\N$
 \[
 M^n_{x}(i)=\sig{j=1}{i}\left[\mathbf{1}_{[V^n_x(t^n_j-)\geq d_j]}-F(V^n_x(t^n_j-))\right].
   \] %Using \eqref{2.1.1}, \eqref{2.3.1}--\eqref{2.6.1}, we obtain the following system equation: for all $t\geq0$, \beq\label{2.7.1}  V^n(t)=($\mathbb A$^n(t)-\lambda^n t)+M^n_v(A(t))-M^n_d(A(t))-\inte{0}{t}F^n(V^n(s-))dA^n(s)+I^n(t),  \eeq where \beq\label{2.8.1}  I^n(t)\equiv \inte{0}{t}\mathbf{1}_{[V^n(s)=0]}(s)ds, \eeq and $I^n(t)$ represents the idle time at the station during time interval $[0,t]$.
Next, define the following centered quantities
\[ S^n(i)\equiv {\frac1n\sig{j=1}{i}(v_j-1)},\quad S^n_{d,x}(i)\equiv \frac1n {\sig{j=1}{i}(v_j-1)}\mathbf{1}_{\{V^n_x(t^n_j-)\geq d_j\}}.
\]
From (\ref{2.1.1}), algebra, and the above definitions, we have for $t\geq0$
\begin{eqnarray}  \label{eq:Vcentered}
V_x^n(t)  & = & \frac{v}{\sqrt{n}} - t  - \frac{1}{\mu^n} \sum_{j=1}^{A^n(t)} F \left( V_x^n(t_j^n-) \right) + \int_0^t \mathbf{1}_{\{ V_x^n(s) = 0 \}} ds  \\
 & & + \frac{n}{\mu^n} \left( \frac{A_x^n(t)}{n} + S^n\left( A_x^n(t) \right) - S_{d,x}^n( A_x^n(t)) - \frac{1}{n}  M_x^n(A_x^n(t)) \right). \nonumber
\end{eqnarray}

%\subsubsection*{Diffusion Scaling} \label{subsubsection:scaled}
%\subsubsection*{Regulator Mapping Representation under Diffusion Scaling} \label{subsubsection:scaled}
%As described in \cite{ReedWard2008}, the basic equation of the offered waiting time process $\{V^n(t):t\geq0\}$ is given by  \beq\label{3.1.1} V^n(t)=\frac{1}{n}\sig{j=1}{A^n(t)}v_j\mathbf{1}_{[V^n(t_j^n-)<d_j^n]} - \inte{0}{t} \mathbf{1}_{[V^n(s)>0]}(s)ds.\eeq
%{\blue(Chihoon, I am updating the writing slightly in the below. One reason is that we should not use our regulator mapping to write $\wt{V}_x^n$ until we know the mapping exists and is unique. Also, I don't think we need to change any of the definitions the way I have set up the elapsed service time, but I am leaving that until we agree.)}
With the initial state $\wt{\mathbb X}^n(0)=(\tau, v)\equiv x\in\mathbb S$, define fluid-scaled and diffusion-scaled quantities to carry out our analysis. For $t\geq0$, let \beq\label{ds}\bar A^n_x(t)\equiv \frac{A^n_x(t)}{n}, \quad \wt A^n_x(t)\equiv \sqrt n\left(\frac1n A^n_x(t)- \lambda (t - \frac{\tau}{\sqrt n}\wedge t) \right), \quad \wt S^n(t) \equiv \sqrt n S^n([nt]), \eeq \[ \wt S^n_{d,x}(t) \equiv \sqrt n S^n_d([nt])\,, \quad \wt M^n_{x}(t) \equiv \frac{1}{\sqrt n} M^n_{x}([nt]). \]

\noindent {Algebra, (\ref{eq:Vcentered}) and substitution of the above scaled quantities into the scaled offered waiting time process
\[
 \widetilde{V}_x^n(\cdot) \equiv \sqrt{n} V_x^n(\cdot),
\]
shows that for $t\geq0$
}
\beq\label{vne}
\widetilde{V}^n_x(t) = v + b^nt + \wt N^n_x(t)  - {\frac{n}{\mu^n}\int_0^t {\sqrt{n}} F \left( \frac{1}{\sqrt{n}} \widetilde{V}^n_{{x}}(s^-) \right)  d \bar{A}^n_x(s)} + \widetilde{I}^n_x(t),
\eeq
where
\begin{eqnarray*}
b^n & \equiv & {\left(\frac{n}{\mu^n}\right)\sqrt n\left(\frac{\lambda^n}{n} -\frac{\mu^n}{n}\right) } \\
\wt N^n_x(t) & \equiv & \left( \frac{n}{\mu^n} \right)\left(\wt S^n(\bar A^n_x(t)) -  \wt S^n_{d,x}(\bar A^n_x(t)) - \wt M^n_{x}(\bar A^n_x(t))+ \wt A^n_x(t) \right),  \\
%%\red{\widetilde{d}^n_x} &\equiv& \left( \frac{n}{\mu^n} \right) \frac{1}{\sqrt n} F\left( \frac{1}{\sqrt{n}} \widetilde{V}^n_{x}(t^n_1-) \right) - \left( \frac{n}{\mu^n} \right) \frac{1}{\sqrt n} \E(v_1) F\left( \frac{1}{\sqrt{n}} \widetilde{V}^n_{x}(t^n_1-) \right)  \\
%\widetilde{\epsilon}^n_x(t) & = &  {\frac{n}{\mu^n}\int_0^t {\sqrt{n}} F \left( \frac{1}{\sqrt{n}} \widetilde{V}^n_{x}(s^-) \right) ds} - \black{\frac{\black{n}}{\mu^n}}\int_0^t {\sqrt{n}} F \left( \frac{1}{\sqrt{n}} \widetilde{V}^n_{\black{x}}(s^-) \right) d \bar{A}^n_x(s), \\
\wt I^n_x(t) & \equiv & \sqrt n \inte{0}{t}\mathbf{1}_{[\wt V^n_x(s)=0]}ds\,.
\end{eqnarray*}

\subsection{Proof of Proposition \ref{moment1}}
%\blue{(Chihoon, I prefer proofs where the authors outline the structure at the beginning, so the reader can know where the proof is going.  So I have updated the ordering. Also, we need to be really careful about the constants and specifying the $N_0$, etc., as that is how we made a mistake before.  Those details still need another read. I am only putting this comment in blue, rather than try to put in blue everywhere something has changed.)}

%For any $p \in (0,q]$,
%Notice the upper bound
%\begin{eqnarray*}
%    \left| \wt{\mathbb X}^n_x (t|x|) \right|^p & = & \left|   \wt\tau_x^n( t|x| ) + \wt V^n_x(t|x|) \right|^p \\
%    & \leq & \left( 1 \vee 2^{p-1} \right) \left(  \wt\tau_x^n( t|x| )^p + \wt V^n_x(t|x|)^p \right)
%\end{eqnarray*}
%follows
{We will establish the claim when $q\in[2,p)$. Then the claim with $q\in[1,2)$ follows from Jensen's inequality.}
From the inequality (cf. Lemma 2 on page 98 in \cite{roussas2014introduction})
\begin{equation*} % \label{eq:CRinequality}
 (a+b)^r \leq (1 \vee 2^{r-1}) (a^r + b^r) \,\,\mbox{ for }\,\, a,b,r \geq 0,
\end{equation*}
%Furthermore, since by definition $\nu^n_x( s ) \in \{0,1\}$ for all $s \geq 0$,
%Taking expectations in the above display shows that
we obtain for $q\in[2,p)$ that
\[
 \E \left[ \left| \wt{\mathbb X}^n_x (t|x| )  \right|^q \right] \leq  2^{q-1} \left(\E \left[ \wt\tau_x^n( t|x| )^q \right] + \E \left[ \wt V^n_x(t|x|)^q \right] \right).
\]
Therefore, it is sufficient to show that there exists  $t_0 \in \R_+$ such that for all $t \geq t_0$,
\begin{equation} \label{eq:to-show-1}
\lim_{|x|\ra\infty}\sup_{n}\frac{1}{|x|^{{q}}}\E\left[\wt\tau^n_x(t|x|)^{{q}} \right]=0
\end{equation}
and
\begin{equation} \label{eq:to-show-2}
\lim_{|x|\ra\infty}\sup_{n} \frac{1}{|x|^{{q}}}\E\left[\wt{V}^n_x(t|x|)^{{q}} \right]=0
\end{equation}
We first show (\ref{eq:to-show-1}) and second show (\ref{eq:to-show-2}).

\begin{proof}[\bf Proof of \eqref{eq:to-show-1}] We begin by observing that by definition $\wt\tau_x^n(t|x|) \leq \sqrt n u_{(A^n_x(t|x|) +1)} / (\lambda n)$, for all $t \geq 0$ and $x \in \mathbb S$.  Next, define the regular (non-delayed) renewal process  $A_2^n(\cdot)$ via
\[
    A_2^n(t) \equiv \max\left\{ i \in \N: \sum_{j=1}^i \frac{u_{j+1}}{\lambda n} \leq t  \right\}, \quad t\geq0,
\]
where the maximum over an empty set is 0, and observe that
\[
    A_x^n(t) \leq A_2^n(t) + 1.
\]
Now, take $t_0=1$.
Then, for $t\geq t_0$,
\[
 \wt\tau_x^n( t|x| )^q \leq  n^{q/2}\sum_{k=2}^{A^n_x(t|x|) +1} \left( \frac{u_k}{\lambda n} \right)^q   \leq n^{q/2} \sum_{k=2}^{A_2^n(t|x|)+2} \left( \frac{u_k}{\lambda n} \right)^q,
\] where the first inequality above uses the fact that $A^n_x(t|x|)\geq1$ because $t_0|x|=|x|\geq \tau\geq u_1=\tau/\sqrt n$.
From Wald's identity,
\[
 \E \left[ \sum_{k=2}^{A^n_2(t|x|) +2} \left( \frac{u_k}{\lambda n} \right)^q \right] = \E \left[ A^n_2(t |x| )+1 \right] \frac{\E [ u_2^q  ]}{\lambda^q n^q}.
\]
Together, the above two displays imply
\begin{equation} \label{eq:intermediate1}
\E \left[ \wt\tau_x^n( t|x| )^q \right] \leq \frac{1}{\lambda^q} \frac{1}{n^{q/2}} \E [u_2^q] \E \left[ A^n_2(t|x|) +1 \right].
\end{equation}
From the elementary renewal theorem, for any $t >0$ and fixed $n$,
\[
 \frac{\E \left[ A^n_2(t|x|) \right]}{\lambda n t|x|} \rightarrow 1 \mbox{ as } |x| \rightarrow \infty,
\]
and so
\begin{equation} \label{eq:intermediate2}
\E \left[ A^n_2(t|x|) \right] \leq 1 + \lambda n t|x|, \mbox{ for all large enough } |x|.
\end{equation}
Substituting (\ref{eq:intermediate2}) into the right-hand-side of (\ref{eq:intermediate1}) shows that for any $t \geq t_0$ and fixed $n$
\begin{equation} \label{eq:final-step}
\E \left[ \wt\tau_x^n( t|x| )^q \right] \leq \frac{1}{\lambda^{q}} \frac{1}{n^{q/2}} \E [u_2^q] (2 + \lambda n t |x| ) \leq \frac{1}{\lambda^q} \frac{1}{n^{q/2}} 2 \E[u_2^q] \left( 1+ \lambda n t|x| \right),
\end{equation}
where the second inequality follows provided $\lambda n t|x| \geq 1$, which is true for large enough $|x|$ and fixed $n$.
%Then, noting $p <q$, H\"{o}lder's %(Chihoon, not sure why you have been calling this Lyapunov?) \chihoon{Holder's ineq. requires $1/p+1/q=1$, $p,q\geq1$} \amy{I used Holder to derive this; see corollary on page 6 in http://www.math.mcgill.ca/dstephens/556/Handouts/Math556-05-Inequalities.pdf.  Holder seems more well-known, but I don't feel strongly.  }
%inequality establishes
%\[
%\E \left[ \wt\tau_x^n( t|x| )^p \right] \leq \left( \E \left[ \tau_x^n( t|x| )^q \right] \right)^{p/q},
%\]
%and so from (\ref{eq:final-step})
%\begin{equation} \label{eq:final-final-step}
% \E \left[ \wt\tau_x^n( t|x| )^p \right]  \leq \frac{\left( 2 \E [u_2^q] \right)^{p/q} }{\lambda^p n^{p/2}} (1+\lambda n t|x|)^{p/q},
%\end{equation}
%Then, for any $t \geq t_0$ and fixed $n$, for all large enough $|x|$.  Finally, %recall $\wt\tau^n\equiv \sqrt n \tau^n$, and
Finally, (\ref{eq:to-show-1}) follows from (\ref{eq:final-step}) %, %noting that $p/q < p/2<p$ {(since $2<q$)}, and
taking, for example, $t_0 = 1$.
\end{proof}

To complete the proof, we must show (\ref{eq:to-show-2}),
which is more involved than (\ref{eq:to-show-1}),
and proceeds following the approach of
Ye and Yao (\cite{YeYao12-interchange}, Lemma 10 and Proposition 11).
First, we establish two versions of pathwise stability results
(Lemmas \ref{lem-10a} and \ref{lem-10b}),
one for any (fixed) $n$-th system and the other for the whole sequence.
With the moment condition ($\mathbb A$1) on the system primitives,
the pathwise stability results are then turned into the moment stability
in Lemma \ref{lem-10cd}, which finally leads to (\ref{eq:to-show-2}).

\begin{lemma} \label{lem-10a}
\textsc{(Stability of $\wt V^n(\cdot)$ for any (fixed) $n$)} Let $\{r_i \}_{i\geq1}$ be a sequence of numbers such that
$r_i \rightarrow \infty$ as $i \rightarrow \infty$
and assume the sequence of initial states $\{x^i \in {\mathbb{S}}\}_{i\geq1}$
satisfies $|x^i| \le r_i $ for all $i$.
Pick any constant $c>1$.
Then, for any fixed $n$,
the following holds (with probability one),
%as $i\rightarrow \infty$,
\begin{eqnarray}
%  \frac{1}{r_i} \wt \tau_{x^i}^n(r_i t) \rightarrow 0
%   \mbox{~~and~~}
 \lim_{i\rightarrow \infty}  \frac{1}{r_i} \wt V_{x^i}^n(r_i t) = 0 ,
  \mbox{~u.o.c.}  \mbox{~~for~~}t\ge \frac{c}{\sqrt{n}} .
   \label{eq-10a-conv}
\end{eqnarray}
\end{lemma}

\begin{lemma} \label{lem-10b}
\textsc{(Stability of $\wt V^n(\cdot)$)}  Let $\{r_n \}$ be a sequence of numbers such that
 $ r_n  \rightarrow \infty$ as $n\rightarrow \infty$
and assume that the sequence of initial states
$\{x^n\in {\mathbb{S}}\}_{n\geq1}$
satisfies $|x^n|\le r_n$.
 Then, for any $\epsilon >0$,
the following holds (with probability one),
\begin{eqnarray}\label{LL3}
 %\frac{1}{r_n} \wt \tau_{x^n}^n(r_n t) \rightarrow 0
  % \mbox{~~and~~}
 \lim_{n\to\infty}\frac{1}{r_n} \wt V_{x^n}^n(r_n t)= 0,
  \mbox{~u.o.c. for ~}  t\ge \epsilon .
\end{eqnarray}
\end{lemma}

\begin{lemma} \label{lem-10cd}
\textsc{(Moment stability)}  Assume $(\mathbb A1)$ and $(\mathbb A2)$. %Suppose for some $p\ge 2$ and $\delta>0$,
%$\E [u_2 + v_2]^{p+\delta} < \infty$.
% Then, the followings also hold.
\\
(a) Letting $\{r_i\}$ and $\{x^i\}$ as in Lemma \ref{lem-10a},
%Then, the following holds for all $n$,
\begin{eqnarray}
 %   \lim_{i\rightarrow \infty}
   %\E  \frac{1}{r_i^p} \left|\wt {\mathbb{X}}_{x^i}^n(r_i t) \right|^p
   %=
   \lim_{i\rightarrow \infty}
     \E  \frac{1}{r_i^q}
     %\left|\wt V_{x^i}^n(r_i t) + \wt \tau_{x^i}^n(r_i t) \right|^p
   \widetilde V_{x^i}^n(r_i t)^q =0 ,
    \mbox{~~for~~}  t\ge \frac{1}{\sqrt{n}}.
    \label{eq-m-stab-i}
\end{eqnarray}
(b) Letting $\{r_n \}$ and $\{x^n\}$ as in Lemma \ref{lem-10b},
%Then, the following holds,
\begin{eqnarray}
  % \lim_{n\rightarrow \infty}
   %\E \frac{1}{r_n^p} \left|\wt {\mathbb{X}}_{x^n}^n(r_n t) \right|^p
   %=
   \lim_{n\rightarrow \infty}
     \E \frac{1}{r_n^q}
     %\left|\wt V_{x^n}^n(r_n t) + \wt \tau_{x^n}^n(r_n t) \right|^p
   \widetilde V_{x^n}^n(r_n t)^q =0 ,
    \mbox{~~for~~} t>0.
     \label{eq-m-stab-n}
\end{eqnarray}
\end{lemma}

{While the proofs of the above three lemmas are provided
in Section \ref{section:LemmaProof},
we provide some intuitions here, for Lemmas \ref{lem-10a} and \ref{lem-10b}
in particular.
% The key to proving the last one (the uniform integrability) already mentioned above.
Consider the (fixed) $n$-th system in Lemma  \ref{lem-10a}.
During the initial period, if it starts with a large initial state,
say, $\widetilde V_{x^i}^n(0) = r_i$,
% (the first component of $= x^i$),
new arrivals will abandon the service with nearly probability one.
On the other hand, the existing workload, under fluid scaling as in \eqref{eq-10a-conv},
drains (i.e., is processed) at the rate $\sqrt{n}$ approximately.
Therefore, the workload $\widetilde V_{x^i}^n (t)$ will reach
the ``normal'' operating state after the initial period
with an order of $r_i / \sqrt{n}$.
The normal operating state,
scaled by $1/r_i$ (where $r_i \rightarrow \infty$), will be approximately zero,
and this is characterized by the convergence in \eqref{eq-10a-conv}.
In Lemma \ref{lem-10b}, the index $n$ approaches infinity.
For large $n$ and hence large $r_n$, consider the $n$-th (scaled) system and suppose it restarts at a time $t$, i.e., $ \frac{1}{r_n} \widetilde V_{x^n}^n (r_n (t +\cdot))$.
Then, the key observation is similar to the above case;
that is, if the initial state (starting at $t$) is bounded by a constant
(say, $ \frac{1}{r_n} \widetilde V_{x^n}^n (r_n t) \le 1$),
it should be approximately zero after a time of $O(1 / \sqrt{n})$.
Indeed by applying Bramson's hydrodynamic approach,
we are able to bound $ \frac{1}{r_n} \widetilde V_{x^n}^n (r_n t)$ for any time $t$ during the (arbitrarily given) period $[0,T]$.
Therefore, our key observation applies for any time $t$, which will establish \eqref{LL3}.
% Given Lemma \ref{lem-10cd}, the proof of \eqref{eq:to-show-2}
% repeats the one for Proposition 11 of \cite{YeYao12-interchange}
% (and thus is omitted).
% \red{(Could you add a bit more specifics, i.e., basic ideas, any modifications to that proof of Proposition 11, if necessary? This is just to enhance the readability.)}
Lastly, Lemma \ref{lem-10cd} plays a pivotal role in establishing the key moment estimate in \eqref{eq:to-show-2}.} Given  Lemma \ref{lem-10cd}, the proof of \eqref{eq:to-show-2}
repeats the one for Proposition 11 of \cite{YeYao12-interchange}
and is provided below for completeness.
%[Alternatively, this proof can be moved to the the end of next section. -HQ]
%}

\begin{proof}[\bf Proof of \eqref{eq:to-show-2}]
Pick any time $t> 0$.
Suppose \eqref{eq:to-show-2} does  not hold;
then, there exists an $\epsilon_0 >0$ and a sequence of initial
states
  $\{x^i \in \mathbb{S} : i=1,2, \ldots \}$
satisfying $\lim_{i\rightarrow \infty} |x^i| = \infty$
such that
\begin{eqnarray}
\sup_{n} \frac{1}{|x^i|^q} \E
     \left|\widetilde V_{x^i}^n(t|x^i| )\right|^q
    > 2 \epsilon_0 .
   \label{eq-33-uc-01}
\end{eqnarray}
Corresponding to each $x^i$, choose an index in the sequence $\{n\}_{n\geq1}$,
denoted by $n_i$, such that
\begin{eqnarray}
 \frac{1}{|x^i|^q} \E
    \left|\widetilde V_{x^i}^{n_i}(t|x^i|)\right|^q
    > \epsilon_0 .
   \label{eq-33-uc-02}
\end{eqnarray}
We claim that $\{ n_i \}_{i\geq1}$ cannot be bounded. Otherwise, at least
an index, say $n'$, repeats in the sequence for infinitely many times;
this contradicts to Lemma \ref{lem-10cd}(a).
{Otherwise,} without loss of generality,
assume $n_i \rightarrow \infty$ as $i\rightarrow \infty$.
Then, the bound in (\ref{eq-33-uc-02}) contradicts to
Lemma \ref{lem-10cd}(b). {We conclude that the aforementioned sequence of initial states  $\{x^i \in \mathbb{S} : i=1,2, \ldots \}$ satisfying (\ref{eq-33-uc-01}) cannot exist, which implies \eqref{eq:to-show-2} holds. }
\end{proof}

\section{Proofs of the Main Results (Theorems~\ref{stability} and~\ref{convstat})} \label{section:TheoremProof}

%\blue{ (I am thinking the following proof goes early on in this Section.)}
\begin{proof}[\bf Proof of Theorem~\ref{stability}] %From Proposition \ref{moment1} with $p=1$, there exist $\delta\equiv t_0>0$ {and $N_0\geq1$} such that {for each $n\geq N_0$} \beq\label{fluidstab}\lim_{|x|\ra\infty}\frac{1}{|x|}\E|\wt{\mathbb X}^n_x(|x|\delta)|  = 0.\eeq
%We begin by noticing that  %Proposition \ref{moment1} continues to hold with $p=1$:
From Proposition \ref{moment1},
there exists $\delta\equiv t_0>0$ such that  \beq\label{fluidstab}\lim_{|x|\ra\infty}\frac{1}{|x|}\E|\wt{\mathbb X}^n_x(|x|\delta)|  = 0.\eeq
%This is because estimates such as \eqref{eq:final-final-step},     \eqref{eq-m-stab-i}, and      \eqref{eq-m-stab-n} continue to hold for $p=1$; see the proof of Lemma \ref{lem-10cd}.
{We require the following lemma, whose proof follows along the same lines of Proposition 4.8 in Bramson \cite{Bramson:2008}. We provide its proof and the notion of petite set in the Appendix for the sake of completeness.
\begin{lemma}\label{petite-set} Assume $(\mathbb A3)$. The set $C=\{x \in \mathbb S: |x| \leq \kappa\}$ is closed petite for every $\kappa>0$. \end{lemma}
\noindent Given Lemma \ref{petite-set}, Theorem 3.1 of \cite{Dai:1995} implies the Markov process $\mathbb X^n(\cdot)$ {and the scaled process $\wt{\mathbb X}^n(\cdot)$ are} positive Harris recurrent and hence the existence of a unique stationary distribution follows.}  %The moment estimate \eqref{fluidstab}, in conjunction with the arguments in the proof of Theorem 3.1 of \cite{Dai:1995}, implies the state-dependent drift criterion of the embedded Markov chain as in Theorem 2.1(ii) of \cite{MeynTweedie94-state-dep}. Proceeding exactly the same as in the proof of Theorem 3.1 of \cite{Dai:1995}, together with the petite set assumption $(\mathbb{A}3)$, the Markov process $\wt{\mathbb X}^n$ is positive Harris recurrent, which implies the existence of a unique stationary distribution.
\end{proof}

Having Proposition~\ref{moment1} and Theorem~\ref{stability}  at hand, the proof of Theorem~\ref{convstat} follows a similar outline to  that of Theorem 3.1 in \cite{BudhLee07}.  (Recall \cite{BudhLee07} establishes the validity of the heavy traffic stationary approximation for a generalized Jackson network without customer abandonment, assuming the inter-arrival and service time distributions have finite polynomial moments, as in $(\mathbb{A}1)$.)  {A global strategy of the proof is as follows.}
First, Proposition \ref{prop2} below establishes uniform (in $n$) estimates on the expected return time of the general Markov process to a compact set. Second, from such return time estimates, moment bounds for the stationary distributions of $\{\wt{\mathbb X}^n\}$, uniform in the scaling parameter $n$, follow readily, yielding tightness of these distributions.
{Third, the distributional convergence in Theorem~\ref{convstat}(a) follows by combining this tightness property with the known
weak convergence results of $\wt V^n$ in (\ref{eq:Vapprox}) %in \cite{ReedWard2008}
and
\cite{WardGlynn2005}. %(see Section~\ref{section:TheoremProof}).
%By suitably strengthening moment conditions on the underlying model primitives, %(that is, by also assuming $(\mathbb A5)$)
Next, we obtain convergence of moments of stationary distributions, i.e., Theorem~\ref{convstat}(b).}
We begin by providing a general statement concerning strong Markov processes.
\begin{proposition}\label{prop2} \black{(Theorem 3.5 of \cite{BudhLee07}, cf. Proposition
5.4 of \cite{Dai:Meyn:1995})} For $n\geq1$, consider a strong Markov process $\{{\mathbb Y}^n_x(t):t\geq0\}$ with initial condition $x$ on a state space $\mathbb T$.  For $\bar\delta\pos$, define the return time to a compact set $C\subset \black{\mathbb T}$ by $\tau^n_C(\bar \delta)\equiv \inf\{t\geq\bar\delta:{\mathbb Y}^n_x(t)\in C\}$. Let $f:\black{\mathbb T}\ra[0,\infty)$ be a measurable map. For $\bar{\delta}\pos$ and a compact set $C\subset \black{\mathbb T}$, define
\[G_n(x)\equiv \E\left[\inte{0}{\ab{\tau}{n}{C}(\bar{\delta})}f({\mathbb Y}^n_x(t))dt\right], \quad x\in \black{\mathbb T}.\] If $\sup_nG_n$ is everywhere finite and
uniformly bounded on $C$, then there exists a constant ${\eta}\pos$, that is independent of $n$, such that for all $n\in\N$, $t \pos$, $x \in\black{\mathbb T}$,
\beq\label{5.11}
\frac{1}{t}\E[G_n(\ab{{\mathbb Y}}{n}{x}(t))] +
\frac{1}{t}\inte{0}{t}\E[f(\ab{{\mathbb Y}}{n}{x}(s))]ds\leq
\frac{1}{t}G_n(x)+{\eta}.\eeq
\end{proposition}

\begin{proof}[\bf Proof of Theorem \ref{convstat}.]   We first prove (a) and then (b).

 \noindent{\bf Part (a):}  Using standard arguments (cf.  \cite{Gamarnik:2006}), it suffices to establish the tightness of the family of stationary distributions $\{ \pi^n:n\geq1 \}$.  \black{Indeed, the tightness implies every subsequence of $\{ \pi^n: n\geq1\}$ admits a convergent subsequence. Denote a typical limit point by $\wt \pi$ and also define the marginal distribution (corresponding to the limiting stationary distribution of $\wt{V}^n$) $\wt\pi_0$ as $\wt\pi_0(A)\equiv \wt\pi(\R_+\times A)$, $A\in \mathcal B(\R_+)$. Then, as in (\ref{eq:Vapprox}),  we see that the process $\wt{V}^n$, with $\wt{\mathbb X}^n(0)$ distributed as $\pi^n$, converges in distribution to $V$ defined in \eqref{3.58} with $V(0)\sim \wt\pi_0$. The stationarity of $\wt V^n$ implies that $\wt\pi_0$ is a stationary distribution for $V$. Since $V$ has a unique stationary distribution, say $\pi$, it must be $\wt \pi_0=\pi$. }

To prove the desired tightness, \black{it suffices to show} that there exists a positive integer $N$ such that for all $n \geq N$
\begin{equation} \label{eq:TheoremToShow}
\inte{\black{\mathbb S}}{}\black{|y|}\pi^n(dy)\leq \tilde c,
\end{equation}
where $\tilde c\pos$ is a constant independent of $n$. The following arguments proceed according to the same outline as the proof of Theorem 3.2 of \cite{BudhLee07}, but with some details that differ.  \black{A key observation from Proposition~\ref{moment1} is that} there exists $\gamma_0 \in (0,\infty)$ such that, for $t_0$  as in that same proposition,
\beq\label{5.1} \sup_{n} \E |\wt{\mathbb X}^n_x(t_0|x|)|^{q}\leq \frac12|x|^{\textcolor{black}{q}}, \mbox{ for all } x \in C^c, \,\, C\equiv \{x\in \mathbb S: |x|\leq \gamma_0\}. \eeq
Next, we apply Proposition \ref{prop2} above with
\[
\mathbb {Y}^n_x=\wt{\mathbb X}^n_x,\,\,\mathbb T=\mathbb S,\,\, \bar{\delta} \equiv t_0 \gamma_0,\,\, f(x) \equiv 1+|x|^{q-1} \mbox{ for } x \in \mathbb S,\,\, q\in[2,p), \,\, C\equiv \{x\in \mathbb S: |x|\leq \gamma_0\}. \] %\black{where $q=2$. (In the proof of part (b), we will take $q\in[2,p)$.)}
Suppose we can show that there exist $N\in\N$ and $\overline{c} \in (0,\infty)$ such that
\begin{equation} \label{eq:to-show}
\sup_{n \geq N} G_n(x) = \sup_{n\geq N}\E\left[\inte{0}{\tau^n_C(\bar\delta)}(1+|\wt{\mathbb X}^n_x(t)|^{{\textcolor{black}{q}-1} })dt \right] \leq \bar c(1+|x|^{\textcolor{black}{q}}), \quad x\in \mathbb S,\end{equation}
so that the conditions of Proposition~\ref{prop2} are satisfied for the family $\{\wt{\mathbb X}_x^n:n\geq N\}$.  Then, for $x \in \mathbb S$ and $\eta\pos$ as in Proposition \ref{prop2},
\[
\Phi_n(x)\equiv \frac1t G_n(x)-\frac1t \E[G_n(\wt{\mathbb X}^n_x(t))] \geq \frac1t \inte{0}{t}\E(f(\wt{\mathbb X}^n_x(s)))ds-\eta,
\] and thus an expectation with respect to the stationary distribution $\pi^n$ has a lower bound,
\beqy
\inte{\mathbb S}{} \Phi_n(x) \pi^n(dx)
&\geq& \inte{\mathbb S}{}\left(\frac1t \inte{0}{t}\E(f(\wt{\mathbb X}^n_x(s)))ds-\eta\right)\pi^n(dx)\nonumber\\
&=& \inte{\mathbb S}{}\left(\frac1t \inte{0}{t}\E(f(\wt{\mathbb X}^n_x(s)))ds\right)\pi^n(dx)-\eta\nonumber\\
&=&  \frac1t \inte{0}{t}\left(\inte{\mathbb S}{} \E(f(\wt{\mathbb X}^n_x(s))) \pi^n(dx)\right) ds -\eta \label{eqn21}\\
&=&  \inte{\mathbb S}{} f(x) \pi^n(dx)  -\eta, \label{eq:needed1}
\eeqy
where \eqref{eqn21} is from Fubini's theorem and \eqref{eq:needed1} follows from the fact that $\pi^n$ is a stationary distribution.
Furthermore,  if $\Phi_n(x)$ is a bounded function in $x$, then from the definitions of the stationary distribution \black{$\pi^n$ and the function  $\Phi_n(x)$}, it is seen that $0 =\inte{\mathbb S}{}\Phi_n(x)\pi^n(dx)$. Otherwise, if $\Phi_n(x)$ is unbounded, then Fatou's lemma implies that (cf. proof of Theorem 3.2 in \cite{BudhLee07} on page 55, also proof of Theorem 5 in \cite{Gamarnik:2006})
\begin{equation} \label{eq:needed3}
0 \geq \inte{\mathbb S}{}\Phi_n(x)\pi^n(dx).
\end{equation}
%Fubini's theorem and the definition of the \black{stationary} distribution yield that
%\begin{equation} \label{eq:needed3}
%\inte{\mathbb S}{} f(x) \pi^n(dx)= \inte{\mathbb S}{} \frac1t \inte{0}{t}\E(f(\wt{\mathbb X}^n_x(s)))ds\pi^n(dx) =   \frac1t \inte{0}{t}\inte{\mathbb S}{} \E(f(\wt{\mathbb X}^n_x(s))) \pi^n(dx) ds.
%\end{equation}
Finally, it follows from (\ref{eq:needed1})--(\ref{eq:needed3}) that
\beq\label{925}
0 \geq \inte{\mathbb S}{} \Phi_n(x) \pi^n(dx) \geq \inte{\mathbb S}{} f(x) \pi^n(dx) - \eta,
\eeq
which establishes the desired uniform moment bound in (\ref{eq:TheoremToShow}).

To complete the proof of part (a), it only remains to show (\ref{eq:to-show}). %which we do for $p \in [2,\infty)$ because that is needed in the proof of part (b) subsequently.
{The following arguments are similar to those leading to Theorem 3.4 in \cite{BudhLee07}. We note that the terms $\overline{\delta}, f(x)$, and the compact set $C$ are chosen exactly the same way as in the cited theorem.  However, minor modifications are necessary because the Markov process $\wt{\mathbb X}^n_x$ in this paper is defined differently from the Markov process $\hat{X}_x^n$ in \cite{BudhLee07}. This entails to checking the following bound, which corresponds to (38) in \cite{BudhLee07}: there exist $N \in \N$ and $c_0 \in (0,\infty)$ such that
\beq\label{tso}
\sup_{n \geq N} \E \left[ \int_0^{\sigma_1} \left( 1 + \wt{\mathbb X}^n_x(t) \right)^{q-1} dt \right] \leq c_0 \left( 1+|x|^q \right),\quad x \in \mathbb S,
\eeq
where $\sigma_1 \equiv t_0 \left( |x| \vee \gamma_0 \right)$ is defined exactly as in the proof of Theorem 3.4 in \cite{BudhLee07}.  Notice that  $\sigma_1\leq c_1(1+|x|)$ for some $c_1\pos$, which is the same estimate as in (39) of \cite{BudhLee07}.}  For  $t\geq0$ and a real-valued function $f$ on $\R_+$, define $\| f \|_t \equiv \sup_{0 \leq s \leq t} |f(s)|$.  To show \eqref{tso}, it is sufficient to show there exist constants $c_2, c_3 \in (0,\infty)$, and $N \in \N$ such that, for all $n \geq N$,
\begin{equation} \label{eq:Vbound}
\E \left \| \widetilde{V}_x^n \right\|^{q-1}_{c_1(1+|x|)} \leq c_2 \left( 1 + |x| \right)^{q-1}
\end{equation}
and
\begin{equation} \label{eq:tauBound}
\E \left\| \wt\tau_x^n \right\|_{c_1(1+|x|)}^{q-1} \leq c_3 (1+|x|)^{q-1}.
\end{equation}

\noindent {\em The argument to show (\ref{eq:Vbound}).} {It follows from the proof of Lemma~\ref{lem-10cd} (more precisely, \eqref{38n} in Section \ref{section:LemmaProof} and the fact that $||\widetilde{V}_x^n||_t\leq ||\widetilde{W}_x^n||_t$, where $\{{W}_x^n(t):t\geq0\}$ denotes the offered waiting time process of the $GI/GI/1$ queue defined as in the proof of Lemma~\ref{lem-10cd}).}

\noindent {\em The argument to show (\ref{eq:tauBound}).}  The same logic used to show (\ref{eq:intermediate1}) also shows
\[
\E \left[ \wt\tau_x^n(t)^q \right] \leq \frac{1}{\lambda^q} \frac{1}{n^{q/2}} \E \left[ u_2^q \right] \left( \E \left[ A_2^n(t) \right]+1  \right), \mbox{ for } t \geq 0,
\]
where $A_2^n$ is the regular (non-delayed) renewal process defined in the second paragraph of the proof of Proposition \ref{moment1}.
Since $A_2^n$ is a non-decreasing process,
\[
 \E \left[ \wt\tau_x^n(t)^q \right] \leq \frac{1}{\lambda^q} \frac{1}{n^{q/2}} \E[u_2^q] \left( \E \left[ A_2^n(c_1 + (1+c_1)|x|) \right] +1 \right), \mbox{ for } 0 \leq t \leq c_1(1+|x|).
\]
From the elementary renewal theorem, there exists $N \in \N$ such that for all $n \geq N$,
\[
\E \left[ A_2^n(c_1 + (1+c_1)|x|) \right] \leq \lambda n \left( c_1 + (1+c_1)|x| +1 \right) = \lambda n (1+c_1)(1+|x|),
\]
and so
\[
 \E \left[\wt\tau_x^n(t)^q \right] \leq \frac{1}{\lambda^q} n^{1-q/2} \E[u_2^q] \left( \lambda(1+c_1)(1+|x|)+1 \right), \mbox{ for } n \geq N \mbox{ and } 0 \leq t \leq c_1(1+|x|).
\]
Recalling that $q \geq 2$ and $n \geq 1$ ensures $n^{1-q/2} \leq 1$ and so
\[
\E \left[\wt \tau_x^n(t)^q \right] \leq \frac{1}{\lambda^q} \E[u_2^q] \left( \lambda(1+c_1)+1 \right)(1+|x|), \mbox{ for } n \geq N \mbox{ and } 0 \leq t \leq c_1(1+|x|).
\]
Since $0 < q-1 <q$, H\"{o}lder's inequality implies
\[
 \E \left[ \wt\tau_x^n(t)^{q-1} \right] \leq \left( \E \left[ \wt\tau_x^n(t)^q \right] \right)^{(q-1)/q}, \mbox{ for } t \geq 0.
\]
Hence for $c_3 \equiv \left( \E[u_2^q] \left( \lambda(1+c_1) +1 \right) / \lambda^q \right)^{(q-1)/q} $,
\[
  \E \left[ \wt\tau_x^n(t)^{q-1} \right] \leq c_3 \left(1+|x| \right)^{(q-1)/q} \leq c_3 (1+|x|)^{q-1}, \mbox{ for } n \geq N \mbox{ and } 0 \leq t \leq c_1(1+|x|),
\]
from which (\ref{eq:tauBound}) follows.

 \noindent{\bf Part (b):} {Proposition~\ref{moment1} implies the moment estimate \eqref{5.1} with $q\in[1,p)$. However, when applying Proposition \ref{prop2} with $f(x) \equiv 1+|x|^{q-1}$, we require $q\in[2,p)$ in order to obtain  (\ref{eq:tauBound}) and \eqref{925}.} Therefore, we have a uniform (in $n$) moment bound on the $(q-1)$-th moment of the family of stationary distributions $\{\pi^n\}_{n\geq1}$.
{If we pick $q$ such that $m<q-1$, this uniform moment bound, in turn, implies} the uniform integrability of $\{|\wt{\mathbb X}^n|^{m}\}_{n\geq1}$ with $m\in(0,q-1)$ in stationarity (i.e., $\wt{\mathbb X}^n(0)\sim \pi^n$). Combining the weak convergence result established in part (a) and the aforementioned uniform integrability, we conclude the desired moment convergence  in stationarity as in part (b). This completes the proof. \end{proof}

\section{Proofs of Lemmas~\ref{lem-10a}--\ref{lem-10cd}} \label{section:LemmaProof}

\begin{proof}[\bf {Proof of Lemma~\ref{lem-10a}.}]
Without loss of generality, assume that as $i \rightarrow \infty$,
$x^i / r_i \rightarrow {\bar x}\equiv (\bar \tau, \bar v(0))$
with $|{\bar x}| = \bar \tau + \bar v(0) \le 1$; otherwise, it suffices
to consider any  convergent subsequence. %(Note the state descriptor $\nu\in\{0,1\}$, an indicator function corresponding to ``empty'' and ``busy'' queue, plays no role in the asymptotic analysis, so we do away with it in this proof.)
Fix the index $n$ throughout the proof.
We also omit the index $n$ and the subscript $x^i$
whenever it does not cause any confusion.

%For any time $t \ge c/\sqrt{n}$, we have,
%as $i\rightarrow \infty$,
%  $$ \frac{r_i t}{ \tau^n(0)}
%   \ge \frac{r_i c/\sqrt{n}}{ \tau^i/\sqrt{n}}
%   \rightarrow \frac{ c}{\bar \tau} \ge c >1 ,
%  $$
%where the first inequality is because
% $\tau^n(0) = \wt \tau_{x^i}^n(0) /\sqrt{n}
%  = \tau^i/\sqrt{n}$.
%%   and the convergence is taken as $i\rightarrow \infty$.
%Hence, the residual time
% $\wt \tau_{x^i}^n(r_i t) [= \sqrt{n} \tau^n(r_i t)] $
%is part of an interarrival time {\it other than}
%the initial residual arrival time.
%Then, the first convergence in (\ref{eq-10a-conv})
%follows from Lemma 5.1 of \cite{Bramson:1998}.

For the $i$-th copy of the $n$-th system,
write the offered waiting time as:
\begin{eqnarray}
 \frac{1}{r_i} \wt V_{x^i}^n(r_i t)  & \equiv &
    \wt v_i(t) = \phi_i(t) + \eta_i(t),
    \mbox{~~~~with}
   \label{eq-lee-ward-8} \\
 \phi_i(t) &  \equiv &
   \frac{v^i}{r_i}
   % + \frac{n\sqrt{n}}{\mu^n} \cdot \frac{A^n(r_i t)}{r_i n}
   + \frac{n\sqrt{n}}{\mu^n} \cdot \frac{1}{r_i n}
     \sum_{j=1}^{A^n(r_i t)} (1- F(V^n(t_j^n -)) )
   - \sqrt{n} t
     \nonumber \\
 && + \frac{n\sqrt{n}}{\mu^n} \cdot \frac{1}{r_i} \left(
       S^n(A^n(r_i t)) - S_{d}^n(A^n(r_i t))
       - \frac{1}{n} M^n(A^n(r_i t))
     \right),
   \label{eq-lee-ward-8b} \\
 \eta_i(t) &  \equiv &
   \frac{\sqrt{n}}{r_i} \int_0^{r_i t} {\bf 1}_{\{V^n(s)=0\}} ds .
   \label{eq-lee-ward-8c}
\end{eqnarray}

First, estimate the item associated with the arrival in the above
(in the first summation):
\begin{eqnarray}
 \frac{A^n(r_i t)}{r_i n}
   &=& \frac{1}{r_i n} \left( A^n(r_i t)
     - \lambda n (r_i t - \frac{\tau_i}{\sqrt{n}} \wedge r_i t )
     \right)
     + \lambda (t - \frac{\tau_i}{r_i \sqrt{n}} \wedge t )
     \nonumber \\
 && \rightarrow
   \lambda (t - \frac{\bar \tau}{\sqrt{n}} \wedge t ),
   \mbox{~~~~as } i \rightarrow \infty
   \mbox{~~a.s.}
     \label{eq-arr-fluid}
\end{eqnarray}

Second, denote the term associated with the arrival and abandonment as
\begin{eqnarray}
  && \xi_i(t)
    \equiv \frac{1}{r_i n} \sum_{j=1}^{A^n(r_i t)} (1- F(V^n(t_j^n -))).
   \label{eq-xibar}
\end{eqnarray}
Observe that for any $0\le t_1 < t_2 $, we have
\begin{eqnarray}
  &&  0 \le \xi_i(t_2) - \xi_i(t_1)
   \le \frac{1}{r_i n} ( A^n(r_i t_2) - A^n(r_i t_1) ).
    \label{eq-xi-lowbdd}
\end{eqnarray}
From (\ref{eq-arr-fluid}), we note that the right-hand side in
the above converges uniformly to
 $\lambda (t_2 - t_1
  - \frac{\bar \tau}{\sqrt{n}} \wedge t_2
  + \frac{\bar \tau}{\sqrt{n}} \wedge t_1 )$.
Therefore, any subsequence of $i$ contains a further subsequence such
that as $i\rightarrow \infty$ along the further subsequence,
we have the weak convergence
\begin{eqnarray}
  && \xi_i(\cdot) \Rightarrow {\bar \xi}(\cdot) , \,\,\mbox{ in }\,\, D(\R) \,\,\mbox{ as } \,\, i\rightarrow \infty,
    \nonumber
\end{eqnarray}
where the limit ${\bar \xi}(\cdot)$ is Lipschitz continuous (recall \eqref{eq-xi-lowbdd})
 with a Lipschitz constant $\lambda$ (with probability one).
Without loss of generality, we can assume the above convergence is along the full sequence,
and furthermore, by using the coupling technique, we can further assume
the convergence is almost surely:
\begin{eqnarray}
  && \xi_i(t) \rightarrow {\bar \xi}(t) ,
    \mbox{~~~~as } i\rightarrow \infty
    \mbox{~~~a.s}.
    \nonumber
\end{eqnarray}

Third, for the martingale terms, we have
as $i\rightarrow \infty$ with probability one,
\begin{eqnarray}
 && \frac{1}{r_i } \left(
       S^n(A^n(r_i t)) - S_{d}^n(A^n(r_i t))
       - \frac{1}{n} M^n(A^n(r_i t)
     \right)
     \rightarrow 0 .
\end{eqnarray}
% \footnote{CH: For the convergence of the $M^n()$ component, I am not completely sure whether a proof is necessary. Convergence of the first two terms (to zero) should be obvious (I think). -HQ}

Putting the above convergences together yields,
as $i\rightarrow \infty$,
\begin{eqnarray}
 && \phi_i(t) \rightarrow  \bar\phi(t)  \equiv
   \bar v(0)
   % + \frac{n\sqrt{n}}{\mu^n} \cdot \lambda (t - \frac{\bar \tau}{\sqrt{n}} \wedge t )
   + \frac{n\sqrt{n}}{\mu^n} {\bar \xi}(t)
   - \sqrt{n} t ,
   \mbox{~~~u.o.c. of } \,\,t\ge0.
     \label{eq-xi-bar-conv}
\end{eqnarray}
% This also implies
% \red{[I think there could be more explanation in how you go from (29) to (30).  After you have (29), then you use the continuity of the integral mapping (or something like this, \blue{[I think the continuity of the reflection term can be found in some reference?]}) to show $\eta_i(t) \rightarrow \overline{\eta}(t)$.  Then you use the expression in (24) to show the convergence to $\overline{v}$.  So it seems like there are two steps hidden there]}
Note from \eqref{eq-lee-ward-8}--\eqref{eq-lee-ward-8c} that
the tuple $( \wt v_i(t), \phi_i(t), \eta_i(t))_{t\geq0}$ satisfies the
one-dimensional linear Skorokhod problem (cf. \S6.2 of \cite{ChenYao01}):
\begin{eqnarray}
 && \wt v_i(t) = \phi_i(t) + \eta_i(t) \ge0, ~~
    d \eta_i(t) \ge0
      \mbox{~with~ } \eta_i(0)=0,  ~~
  \wt v_i(t) d \eta_i(t) = 0 .
    \nonumber
\end{eqnarray}
% \begin{eqnarray}
%  && \wt v_i(t) = \phi_i(t) + \eta_i(t) \ge0, ~~
%     \nonumber \\
%  &&  d \eta_i(t) \ge0, ~ \eta_i(0)=0
%     \nonumber \\
%  && \bar v_i(t) d \bar\eta_i(t) = 0 .
%     \nonumber
% \end{eqnarray}
Hence, by invoking the Lipschitz continuity of the Skorokhod mapping
(cf. Theorem 6.1 of \cite{ChenYao01}),
the convergence in (\ref{eq-xi-bar-conv}) implies
\begin{eqnarray}
 && \frac{1}{r_i} \wt V_{x^i}^n(r_i t)
   \rightarrow \bar v(t)
   \mbox{~~and~~}
   \eta_i(t) \rightarrow  \bar\eta(t)
   \mbox{~~~u.o.c. of }\,\, t\ge0,
     \label{eq-V-fluidlimit}
\end{eqnarray}
with the limit satisfying the Skorokhod problem as well:
\begin{eqnarray}
 && \bar v(t) = \bar\phi(t) + \bar\eta(t) \ge0, ~~
   d \bar\eta(t) \ge0
     \mbox{~with~} \bar\eta(0) =0, ~~
   \bar v(t) d \bar\eta(t) = 0 .
     \label{eq-sko_fluid}
\end{eqnarray}
%\red{[Also, to get to  \eqref{eq-sko_fluid}, don't you need to know the Skorokhod problem has a unique solution?  That probably should be said.
 % \blue{[The uniqueness of the solution is not used here. -HQ]}
%]}

Next, we further examine the limit $\bar \xi(\cdot)$
following the approach of Chen and Ye
(\cite{ChenYe12}, Proposition 3(b)).
From (\ref{eq-arr-fluid}) and (\ref{eq-xibar}),
and noting that
 $ \xi_i(t) \le \frac{A^n(r_i t)}{r_i n}$,
we have
\begin{eqnarray}
  \bar \xi(t) = 0, ~~ 0\le t \le \frac{\bar \tau}{\sqrt{n}} .
   \label{eq-xi-inittime}
\end{eqnarray}
Now, consider any regular time $t_1 > {\bar \tau}/{\sqrt{n}}$,
at which all processes concerned, i.e.,
$\bar v(\cdot), \bar\phi(\cdot), $ and $\bar\eta(\cdot)$ are differentiable,
and $\bar v(t_1) >0$.
Note that the Lipschitz continuity of $\bar \xi_i(\cdot)$
implies that $(\bar v(\cdot), \bar\phi(\cdot), \bar\eta(\cdot))$ are also
Lipschitz continuous.
Therefore, we can find (small) constants $\epsilon>0$ and $\delta>0$
such that the following inequality holds for all sufficiently large
$i$:
\begin{eqnarray}
 \wt v_i(t_2) > \epsilon
   \mbox{~~i.e., }
   \wt V^n(r_i t_2) > r_i \epsilon ,
   ~~~~ t_2 \in [t_1, t_1 + \delta) .
    \label{eq-vi-epsilon}
\end{eqnarray}
Observe that if the $j$-th arrival falls between
$A^n(r_i t_1)+1$ and  $A^n(r_i t_2)$,
then its arrival time, $t_j^n$, shall also falls between the
corresponding time epochs, i.e.,
$r_i t_1 < t_j^n \le r_i t_2$.
Given the estimate in (\ref{eq-vi-epsilon}),
this implies the following estimate holds:
%for those times $t_j^n$:
  $$\wt V^n(t_j^n) > r_i \epsilon .$$
Consequently, we have for all sufficiently large $i$ that
%and all $t_2 \in [t_1, t_1 + \delta)$:
\begin{eqnarray}
   \xi_i(t_2) - \xi_i(t_1)
    &=& \frac{1}{r_i n} \sum_{j= A^n(r_i t_1)+1}^{A^n(r_i t_2)}
    (1 - F(V^n(t_j^n -)) )
     \nonumber \\
    &\le&  \frac{1}{r_i n} (A^n(r_i t_2)- A^n(r_i t_1) )
   ( 1 - F(r_i \epsilon) ),
   ~~~ t_2 \in [t_1, t_1 + \delta) .
    \nonumber
\end{eqnarray}
Taking $i\rightarrow \infty$, the above yields
  $ \bar \xi(t_2) - \bar \xi(t_1) \le 0, $
which is effectively $ \bar \xi(t_2) - \bar \xi(t_1) = 0 .$
% On the other hand, by reviewing the inequality in (\ref{eq-xi-lowbdd}),
% we can show that
%   $ \bar \xi(t_2) - \bar \xi(t_1) \le \lambda(t_2 - t_1), $
% Therefore, we have for all $t_2 \in [t_1, t_1 + \delta)$,
%   $ \bar \xi(t_2) - \bar \xi(t_1) = \lambda(t_2 - t_1) $.
In summary, the above implies
for any regular time $t > {\bar \tau}/{\sqrt{n}}$
with $\bar v(t) >0$,
\begin{eqnarray}
 \frac{d {\bar \xi}(t)}{dt} = 0 .
   \label{eq-xi-deri}
\end{eqnarray}

Finally, using the properties in
(\ref{eq-xi-bar-conv}), (\ref{eq-sko_fluid}), (\ref{eq-xi-inittime})
and (\ref{eq-xi-deri}),
it is direct to show that
 $\bar v(t) = \bar v(0) - \sqrt{n}t $
for $t \le \bar v(0)/ \sqrt{n} (\le 1 / \sqrt{n})$
and
 $\bar v(t) = 0 $
for $t \ge \bar v(0)/ \sqrt{n}$.
This property, along with the convergence in (\ref{eq-V-fluidlimit})
yields the desired convergence in (\ref{eq-10a-conv}). % \hfill $\square$ %$\qed$
\end{proof}

\begin{proof}[\bf {Proof of Lemma~\ref{lem-10b}.}]

We adopt Bramson's hydrodynamics approach
(cf. \cite{Bramson:1998,ms,YeYao-multi-bn}),
and its variation (in Appendix B.2 of \cite{YeYao12-interchange})
in particular,
to examine the processes involved in  \eqref{LL3}. % Lemma \ref{lem-10b}.
Define for each $n$ and for $j=0, 1, \ldots$,
\begin{eqnarray}
 \bar V^{n,j}(u)
  = \frac{1}{r_n} \widetilde V_{x^n}^n(r_n \frac{j+u}{\sqrt{n}} ),
%    = \frac{\sqrt{n}}{r_n} V^n(r_n \frac{j+u}{\sqrt{n}})
%    \label{eq-hydro-1} \\
% \mbox{~~and~~}
% \bar \tau^{n,j}(u)
%  = \frac{1}{r_n} \wt \tau_{x^n}^n(r_n \frac{j+u}{\sqrt{n}} ) ,
%%    = \frac{\sqrt{n}}{r_n} \tau^n(r_n \frac{j+u}{\sqrt{n}})
   ~~~~ u\ge0.
     \label{eq-hydro-2}
\end{eqnarray}
That is, the $n$-th diffusion-scaled process
(say, $\{\frac{1}{r_n} \wt V_{x^n}^n(r_n t )\}_{t\geq0}$)
breaks into many pieces of fluid-scaled process
(i.e., $\{\bar V^{n,j}(u): u\in[0,1]\}$ and $j=0,1, \ldots$),
with each piece covering a period of $1/\sqrt{n}$ in the
diffusion-scaled  process.
(As in the proof of the previous Lemma~\ref{lem-10a}, we omit the index $n$
and the subscript $x^n$ whenever it does not cause any confusion.) %\chihoon{Which index $n$? Did you mean a superscript $n$? But, it shows many times below?}

Let $j_n$ be any nonnegative integer for each $n$
and consider any subsequence of positive integer, denoted ${\cal N}$.
If
 $\limsup_{n\rightarrow \infty, n\in {\cal N}}
   [\bar V^{n,j_n}(0) %+ \bar \tau^{n,j_n}(0)
   ] \le 1 ,
 $ %
% then for any subsequence of integer, there exists a further
% subsequence, denoted ${\cal N}$, such that
then it can be seen that: as $n\rightarrow \infty$ along ${\cal N}$,
\beq\label{L5}
 \bar V^{n,j_n}(u) %+ \bar \tau^{n,j_n}(u)
   \rightarrow 0,
   \mbox{~~~~u.o.c. of }
   u \ge 1 .
\eeq
This can be proved in the same manner as Lemma \ref{lem-10a},
with extra modification on the shifted initial residual arrival time
and the sequence of scaling factors
($\{ r_n / \sqrt{n} \}$ here versus $\{r_i\}$ of Lemma \ref{lem-10a}).
% The detailed proof is omitted.
We defer the proof of \eqref{L5} for now, and proceed to show the main claim in \eqref{LL3}. %We outline a proof later.}
% \red{(Again, could you add a bit more specifics here, i.e., basic ideas, any modifications to that proof of Lemma \ref{lem-10a}, as necessary?)}

%\red{\begin{lemma} \label{lem-10b-hydro}
%Let $j_n$ be any nonnegative integer for each $n$;
%and consider any subsequence of integer, denoted ${\cal N}$.
%If
% $$\limsup_{n\rightarrow \infty, n\in {\cal N}}
%   [\bar V^{n,j_n}(0) %+ \bar \tau^{n,j_n}(0)
%   ] \le 1 ,
%   $$ %
%% then for any subsequence of integer, there exists a further
%% subsequence, denoted ${\cal N}$, such that
%then the following convergence holds
%as $n\rightarrow \infty$ along ${\cal N}$,
%\begin{eqnarray}
% \bar V^{n,j_n}(u) %+ \bar \tau^{n,j_n}(u)
%   \rightarrow 0,
%   \mbox{~~~~u.o.c. of }
%   u \ge 1 .
%     \label{eq-hd-conv-lem}
%\end{eqnarray}
%\end{lemma}}
%
%
%\begin{proof}[\bf {Proof of Lemma~\ref{lem-10b-hydro}.}]
%The lemma is proved in the same manner as Lemma \ref{lem-10a},
%with extra care on the shifted initial residual arrival time
%and the sequence of scaling factors
%($\{ r_n \sqrt{n} \}$ here versus $\{r_i\}$ of Lemma \ref{lem-10a}).
%The detailed proof is omitted.
%\end{proof}
%

To prove \eqref{LL3}, %Lemma \ref{lem-10b},
given the hydrodynamic representation of the waiting time in \eqref{eq-hydro-2},
it suffices to show that
% is implied by the following claim:
for any $\Delta>0$ and $\epsilon>0$,
the following holds for sufficiently large $n$
and for $j= 1, \ldots, \lfloor\sqrt{n} \Delta\rfloor$ (excluding $j=0$),
\begin{eqnarray}
 \bar V^{n,j}(u) %+ \bar \tau^{n,j}(u)
   \le \epsilon,
   \mbox{~~~~for }
   u\in [0,1] .
     \label{eq-hd-equiv}
\end{eqnarray}
% \red{(Why is the above sufficient to prove Lemma \ref{lem-10b}?)}
We prove this by contradiction.
Suppose to the contrary, there exists a subsequence ${\cal N}$,
and for all $n\in {\cal N}$, we can find an integer index
$j_n \in [1, \sqrt{n} \Delta]$ and a time $u_n\in [0,1]$
such that
\begin{eqnarray}
 \bar V^{n,j_n}(u_n) %+ \bar \tau^{n,j_n}(u_n)
 > \epsilon .
   \label{eq-contradict}
\end{eqnarray}
Furthermore, we can require that for each $n\in {\cal N}$,
$j_n$ is the smallest integer
(in $[1, \sqrt{n} \Delta]$) for the above inequality to hold.

Observe that the following initial condition must hold
for sufficiently large $n\in {\cal N}$:
\begin{eqnarray}
 \bar V^{n,j_n-1}(0)% + \bar \tau^{n,j_n-1}(0)
 \le 1 .
     \label{eq-hd-c0}
\end{eqnarray}
First, if $j_n=1$, then from the definition in
% (\ref{eq-hydro-1}) and
(\ref{eq-hydro-2}), we evaluate the left-hand side
of the above as:
\begin{eqnarray}
 \bar V^{n,0}(0) %+ \bar \tau^{n,0}(0)
  = \frac{1}{r_n} \widetilde V_{x^n}^n(0) %+ \wt \tau_{x^n}^n(0))
  \le \frac{|x^n|}{r_n}
  \le 1.
     \nonumber
\end{eqnarray}
On the other hand, we have $j_n \ge 2$;
and according to the definition of $j_n$,
it means that the inequality (\ref{eq-hd-equiv}) applies
for $j= j_n-1$. As a result, the inequality (\ref{eq-hd-c0})
also holds.

% \frac{1}{r_n} \wt V_{x^n}^n(r_n \frac{j+u}{\sqrt{n}} )
%    = \frac{\sqrt{n}}{r_n} V^n(r_n \frac{j+u}{\sqrt{n}})
%    \label{eq-hydro-1} \\
%  \bar \tau^{n,j}(u)
%   & \equiv& \frac{1}{r_n} \wt \tau_{x^n}^n(r_n \frac{j+u}{\sqrt{n}} )
%    = \frac{\sqrt{n}}{r_n} \tau^n(r_n \frac{j+u}{\sqrt{n}})
%      \label{eq-hyd

Given the above initial condition,
applying \eqref{L5} %Lemma \ref{lem-10b-hydro}
yields the following convergence
as $n\rightarrow \infty$ in ${\cal N}$,
\begin{eqnarray}
 \bar V^{n,j_n-1}(u) %+ \bar \tau^{n,j_n-1}(u)
   \rightarrow 0,
   \mbox{~~~~u.o.c. of } u\ge 1.
     \nonumber
\end{eqnarray}
Therefore, for sufficiently large $n\in {\cal N}$, we have
\begin{eqnarray}
 \bar V^{n,j_n}(u) %+ \bar \tau^{n,j_n}(u)
   = \bar V^{n,j_n-1}(1+u) %+ \bar \tau^{n,j_n-1}(1+u)
    \le \epsilon,
    ~~~~ u \in [0,1],
     \nonumber
\end{eqnarray}
which contradicts to (\ref{eq-contradict}), and this establishes the desired result in \eqref{LL3}.

Finally, it remains to show the claim  \eqref{L5} holds,
which resembles the proof of Lemma \ref{lem-10a} and hence we provide the outline only.
%\noindent {\bf Proof (outline) of \eqref{L5}.}
Without loss of generality, assume that as $n \rightarrow \infty$,
%   $x^n  \equiv \bar V^{n,j_n}(0)
%     \rightarrow {\bar x}\equiv (\bar \tau, \bar v(0))$
% with $|{\bar x}| = \bar \tau + \bar v(0) \le 1$.
  $\bar V^{n,j_n}(0) \rightarrow  \bar v(0) \le 1 .$
Similar to (\ref{eq-lee-ward-8}), the offered waiting time for the $n$-system
can be expressed as
\begin{eqnarray}
 \bar V^{n,j_n}(u)  & = &
   \frac{1}{r_n} \wt V^n(r_n \frac{j_n+u}{\sqrt{n}})
     \equiv \wt v_n(u) \equiv \phi_n(u) + \eta_n(u),
   \label{eq-lee-ward-8b}
\end{eqnarray}   where
\begin{eqnarray}
 \phi_n(u) & \equiv& \bar V^{n,j_n}(0)
   + \frac{n}{\mu^n} \cdot \frac{1}{r_n \sqrt{n}}
     \sum_{j=A^n(r_n j_n /\sqrt{n})+1}^{A^n(r_n (j_n+u)/\sqrt{n})}
        (1- F(V^n(t_j^n -)) )
   - u
     \nonumber \\
 && + \frac{n}{\mu^n} \cdot \frac{\sqrt{n}}{r_n} \left(
       S^n(A^n(r_n \frac{j_n+u}{\sqrt{n}}))
     - S^n(A^n(r_n \frac{j_n}{\sqrt{n}}))
     \right)
     \nonumber \\
 && - \frac{n}{\mu^n} \cdot \frac{\sqrt{n}}{r_n} \left(
         S_{d}^n(A^n(r_n \frac{j_n+u}{\sqrt{n}}))
       - S_{d}^n(A^n(r_n \frac{j_n}{\sqrt{n}}))
     \right)
     \nonumber \\
 && - \frac{n}{\mu^n} \cdot \frac{1}{r_n \sqrt{n}} \left(
         M^n(A^n(r_n \frac{j_n+u}{\sqrt{n}})
       - M^n(A^n(r_n \frac{j_n}{\sqrt{n}})
     \right), % \,\,\mbox{ and }
     \nonumber \\
 \eta_n(u) & \equiv &
   \frac{\sqrt{n}}{r_n}
     \int_{\frac{j_n}{\sqrt{n}}}^{\frac{j_n+u}{\sqrt{n}}}
       {\bf 1}_{\{V^n(s)=0\}} ds .
   \nonumber
\end{eqnarray}

In parallel to the convergence in (\ref{eq-arr-fluid}),
we have:
\begin{eqnarray}
 && \frac{1}{r_n \sqrt{n}}
   \left( A^n(r_n \frac{j_n+u}{\sqrt{n}})
   - A^n(r_n \frac{j_n}{\sqrt{n}}) \right)
    \nonumber \\
 &=& \frac{1}{r_n \sqrt{n}} \left(
   A^n(r_n \frac{j_n+u}{\sqrt{n}}) - A^n(r_n \frac{j_n}{\sqrt{n}})
     - \lambda r_n \sqrt{n} \left( u
        - u \wedge \left( \frac{1}{r_n} {\wt \tau_n}(r_n \frac{j_n}{\sqrt{n}})  \right) \right)
     \right)
     \nonumber \\
 &&   + \lambda \left( u
        - u \wedge \left( \frac{1}{r_n} {\wt \tau_n}(r_n \frac{j_n}{\sqrt{n}})  \right) \right)
     \nonumber \\
 & \rightarrow &
   \lambda u ,
   \mbox{~~~~as } n \rightarrow \infty
   \mbox{~~~a.s}.
     \label{eq-arr-fluid-b}
\end{eqnarray}
While the convergence in (\ref{eq-arr-fluid}) is justified
by the conventional functional strong law of large numbers
for a renewal process (e.g., \cite{ChenYao01}),
we have applied here a version established by
Bramson (\cite{Bramson:1998}, Proposition 4.2)
and Stolyar (\cite{Stolyar04}, Appendix A.2), which accompanies
the hydrodynamic scaling approach.
% Moreover, we have applied the convergence,
%  $ \frac{1}{r_n} {\wt \tau_n}(r_n \frac{j_n}{\sqrt{n}}) \rightarrow 0$,
% in the above, which is due to (\ref{eq:to-show-1}).

Then, following the arguments from \eqref{eq-xibar}--\eqref{eq-sko_fluid},
it is seen that, as $n\rightarrow \infty$,
\begin{eqnarray}
 &&  (\bar V^{n,j_n}(u) , \phi_n(u) , \eta_n(u) )
   \rightarrow (\bar v(u), \bar\phi(u) , \bar\eta(u) )
   \mbox{~~u.o.c. of } u\ge0,
     \label{eq-V-fluidlimit-b}
\end{eqnarray}
where
\begin{eqnarray}
 && \bar\phi(u) =  \bar v(0)
   + \frac{1}{\lambda} {\bar \xi}(u) - u .
     \label{eq-xi-bar-conv-b}
\end{eqnarray}
Moreover, the limits, $\bar \xi(\cdot)$, $\bar v(\cdot)$, $\bar\phi(\cdot)$ and $\bar\eta(\cdot)$,
are Lipschitz continuous and satisfy the Skorokhod problem:
\begin{eqnarray}
 && \bar v(u) = \bar\phi(u) + \bar\eta(u) \ge0, ~~
    d \bar\eta(u) \ge0
      \mbox{~with~ } \bar\eta(0)=0,  ~~
  \bar v(u) d \bar\eta(u) = 0 .    \label{eq-sko_fluid-b}
\end{eqnarray}

Next, a similar analysis on the limit ${\bar \xi}(\cdot)$ yields
% \begin{eqnarray}
%   \bar \xi(u) = 0, ~~ 0\le u \le \bar \tau ,
%    \label{eq-xi-inittime-b}
% \end{eqnarray}
% and,
that, for any regular time $u > \bar \tau$
with $\bar v(u) >0$,
\begin{eqnarray}
 \frac{d {\bar \xi}(u)}{du} = 0 .
   \label{eq-xi-deri-b}
\end{eqnarray}

Finally, using the properties in
(\ref{eq-xi-bar-conv-b}), (\ref{eq-sko_fluid-b})
% (\ref{eq-xi-inittime-b})
and (\ref{eq-xi-deri-b}),
it is direct to show that
 $\bar v(u) = \bar v(0) - u $
for $u \le \bar v(0) (\le 1 )$
and
 $\bar v(u) = 0 $
for $u \ge \bar v(0)$.
This property, along with the convergence in (\ref{eq-V-fluidlimit-b})
yields the desired convergence in the claim (\ref{L5}).
\end{proof}

%\red{[Chihoon: Please relocate the above proof. Two suggestions.
%(1) Postpone it to the end of the proof of Lemma \ref{lem-10b}, or (2) Create another lemma for (\ref{L5}) (like my first draft).]}

\begin{proof}[\bf {Proof of Lemma~\ref{lem-10cd}.}]
{%Letting $q\equiv p+\delta \in(2,\infty)$ from $(\mathbb{A}1)$,
We first show that
 for any (fixed) $p' \in (q, p)$,}
there exists a constant $\kappa'\pos$
such that the following bound holds
for any initial state $x$ and any time $t\ge0$,
\begin{eqnarray}
% \E [\wt {\mathbb{X}}_{x}^n( t)]^{p'}
 \E [\wt {V}_{x}^n( t)]^{p'}
   \le \kappa' (1+ |x|^{p'} + t^{p'}) .
    \label{eq-moment-bdd-X}
\end{eqnarray}
Indeed, if we discard the abandonment component, the systems will
be reduced to the more basic $GI/GI/1$ queues. Dropping the abandonment component in the original state  equation \eqref{2.1.1}, one gets %Consider a sequence of systems that are same as the ones under studied
%but do not abandon jobs.
the offered waiting time process $\{W^n(t):t\geq0\}$, for the $n$-th system, as
\begin{eqnarray}
 W^n_x(t) & = & V^n_x(0)
     + \sum_{j=1}^{A^n_x(t)} v_j^n
     - \int_0^t {\bf 1}_{[W^n_x(s)>0]} ds .
      \nonumber
\end{eqnarray}
Then, {because the system without abandonment dominates the system with abandonment at all times (as can be seen using the one-dimensional Skorokhod mapping),} %Clearly, these systems dominate the original ones (that abandon jobs):
$V^n_x(t) \le W^n_x(t) $ for all $t\ge 0$.
Rewrite the above, with diffusion scaling, as,
\begin{eqnarray}
 && \wt W^n_x(t)  [ \equiv  \sqrt{n} W^n_x(t)]
   = \wt V^n_x(0) + \xi_n(t)
     + \sqrt{n} \int_0^t {\bf 1}_{[\wt W^n_x(s)=0]} ds ,
      \nonumber
\end{eqnarray}
where
\begin{eqnarray}
 \xi_n(t) &\equiv&
  \sqrt{n} \sum_{j=1}^{A^n_x(t)} v_j^n - \sqrt{n} t
     \nonumber \\
 & = & \frac{\sqrt{n}}{\mu^n} \sum_{j=1}^{A^n_x(t)} \left( v_j - 1 \right)
     + \frac{\sqrt{n}}{\mu^n} (A^n_x(t) - n\lambda t)
     + \sqrt{n} \left(\lambda -\frac{\mu^n}{n}\right) \frac{n}{\mu^n}t .
    \nonumber \\
 & = & \frac{n}{\mu^n} ( \wt S^n(\bar A^n_x(t))
    + \wt A^n_x(t) + \theta^n t ) \label{98n}.
\end{eqnarray}
Observe that $\wt W^n_x(\cdot) = \Phi(\wt V^n_x(0) + \xi_n(\cdot))$,
where $\Phi$ is the standard
one-dimensional Skorokhod mapping (cf. \cite{ChenYao01}),
Therefore, it is bounded by the ``free process''
as,
\beq\label{38s}
 \sup_{0\le s\le t}| \wt W^n_x(s) | \le
    | \wt V^n_x(0) |  + 2 \sup_{0\le s \le t} |\xi_n(s)| ,
\eeq
which, under the assumed moment condition, % (in \E Thm 2(b) of lee-ward 2018-04-17),
implies the following bound (whose proof is presented after the next paragraph) for some constant $\kappa'\pos$,
\beq\label{38n}
 \E [\wt V^n_x( t)]^{p'}
   \le \E [\wt W^n_x( t)]^{p'}
   \le  \E[\sup_{0\le s\le t}| \wt W^n_x(s) |^{p'} ]
   \le \kappa' (1+ |x|^{p'} + t^{p'}),
\eeq
%Combining the above with the bound for residual arrival process
%in (48) [of lee-ward 2018-04-17]
%yields the inequality in (\ref{eq-moment-bdd-X}).
which implies (\ref{eq-moment-bdd-X}).

For the conclusion in (a), the above implies that
 the following bound holds uniformly over all (large) $i$,
 $$ \E \frac{1}{r_i^{p'}}
    \left|\wt {{V}}_{x^i}^n(r_i t) \right|^{p'}
   \le \kappa' (2+  t^{p'}) .
 $$
This implies that the sequence
 $\{ \frac{1}{r_i^q} |\wt {V}_{x^i}^n(r_i t) |^q,
  i\in \N \} $
(where $q<p'$)
is uniformly integrable.
Thus, the limit and the expectation in (\ref{eq-m-stab-i})
can be interchanged, and then the conclusion (a) follows from
Lemma \ref{lem-10a}.
The conclusion (b) is proved in the same way by using Lemma \ref{lem-10b}.

Finally, it remains to show the third inequality in \eqref{38n}. For  $t\geq0$ and a real-valued function $f$ on $\R_+$, define $\| f \|_t \equiv \sup_{0 \leq s \leq t} |f(s)|$.  Recalling the definition in \eqref{98n}, it suffices to prove the following bounds \eqref{eq:KTbound} and \eqref{eq:KTbound2}:  %Notice that
\begin{equation} \label{eq:KTbound}
\E \left[ \| \wt{A}^n_x \|_t^{p'} \right] \leq \E \left[ \| \wt{A}^n_2 \|_{t+|x|}^{p'} \right] \leq c_0 \left( 1+\sqrt{t+|x|} \right)^{p'}
\end{equation}
and
\begin{equation} \label{eq:KTbound2}
\E \left[ \| \wt{S}^n \circ \bar A^n_x \|_t^{p'} \right] \leq \E \left[ \| \wt{S}^n \circ \bar A^n_2 \|_{t+|x|}^{p'} \right] \leq c_1 \left( 1+\sqrt{(t+|x|)} \right)^{p'}
\end{equation}
for some constants $c_0, c_1\pos$, independent of $n$. The first inequalities of \eqref{eq:KTbound} and \eqref{eq:KTbound2} follow from recalling
the regular (non-delayed) renewal process $A_2^n(\cdot)$ defined in the second paragraph of the proof of Proposition~\ref{moment1}, for which $A_x^n(t) \leq A_2^n(t) +1$,
%Also recall that
%\[
%\wt N^n_x(t) \equiv \left( \frac{n}{\mu^n} \right)\left( \wt A^n_x(t) + \wt S^n(\bar A^n_x(t)) -  \wt S^n_{d,x}(\bar A^n_x(t)) - \wt M^n_{x}(\bar A^n_x(t))\right).
%\]
and hence, %Assume the fluid and diffusion scaled processes for $A_2^n$, $\bar{A}_2^n$ and $\widetilde{A}_2^n$, are defined as for $A_x^n$, in the paragraph before the GRM Definition \ref{grm}, so that
\[
\bar{A}_x^n(t) \leq \bar{A}_2^n(t) + \frac{1}{n} \mbox{ and } \wt{A}_x^n(t) \leq \wt{A}_2^n(t) + \frac{1}{\sqrt{n}}.
\]
%\red{(Amy:  The below needs updated, in response to this new bound.  I did not updated, because I wanted you to double-check that you agree we need to do this.  Given that we are updated, I observe that $T$ is not defined.  Maybe we should clarify between $t$ and $T$, and decide whether or not we want both.)}

Then the  second inequality in \eqref{eq:KTbound} follows from Theorem 4 (A1.3) in~\cite{KrichaginaTaksar1992}. For the second inequality in \eqref{eq:KTbound2}, from (A1.16) in~\cite{KrichaginaTaksar1992},
\beq\label{33est}
\E \left[ \| \wt{S}^n \circ \bar A^n_0 \|_t^{p'} \right] \leq c_2 \E[\bar A^n_0(t)^{p'/2}],
\eeq
where $c_2\pos$ is a constant independent of $n$.
Since $p'/2<p'$, from Lyapunov's inequality,
\beq\label{34est}
\E[\bar A^n_0(t)^{p'/2}]^{2/p'} \leq \E[\bar A^n_0(t)^{p'}]^{1/p'}\leq c_3(1+t),
\eeq
for some $c_3\pos$, where the second inequality follows from Theorem  4 (A1.1) in~\cite{KrichaginaTaksar1992}.
Thus, we conclude $\E \left[ \| \wt{S}^n \circ \bar A^n_0 \|_t^{p'} \right]^{2/p'} \leq c_2^{2/p'} c_3(1+t)$, which implies \eqref{eq:KTbound2}.
\end{proof}
% \hfill $\square$
%$\qed$

\section*{Appendix: Proofs of Proposition \ref{pc} and Lemma \ref{petite-set}}
\begin{proof}[\bf Proof of Proposition \ref{pc}]
%Theorem 1(a) of  \cite{WardGlynn2005} adapted to the setting in this paper
%\subsection{Proof of Proposition \ref{pc}}
%\begin{remark}\label{ds1}
The weak convergence in \eqref{eq:Vapprox} can be seen from the expression \eqref{vne} and Theorem 1(a) of  \cite{WardGlynn2005} with minor modifications %(\ref{rep0})
as follows.  Observe that
\begin{equation} \label{eq:WeakConvN}
\wt{N}_x^n \Rightarrow \sqrt{ \frac{\mbox{var}(u_2) + \mbox{var}(v_2) }{\lambda} } W, \,\,\mbox{ in }\,\, D(\R)\,\,\mbox{ as }\,\, n \rightarrow \infty,
\end{equation}
where $W$ is a one-dimensional standard Brownian motion. This is because
standard results on renewal processes show that the delayed renewal process defined in \eqref{254} has the limiting behavior \[
    \bar{A}_x^n \rightarrow \lambda e\,\, \mbox{ and }\,\, \wt{A}_x^n \Rightarrow \sqrt{\lambda} \sqrt{\mbox{var}(u_2)} W_1, \,\,\mbox{ in }\,\, D(\R) \,\,\mbox{ as }\,\, n \rightarrow \infty,
    \]
    where $W_1$ is a standard Brownian motion and $e(\cdot)$ represents the identity map, i.e., $e(t) = t$ for all $t \in \R_+$. Donsker's theorem and a random time change show
\[
\wt{S}^n \circ \bar{A}_x^n \Rightarrow \sqrt{\lambda} \sqrt{\mbox{var}(v_2)} W_2, \,\,\mbox{ in }\,\, D(\R) \,\,\mbox{ as }\,\, n \rightarrow \infty,
\]
where $W_2$ is a standard Brownian motion independent of $W_1$.
 Arguments similar to those in Theorem 5.1(a) in \cite{WardGlynn2005} show
\[
    \wt{S}_{d,x} \circ \bar{A}_x^n \Rightarrow 0 \,\,\mbox{ and }\,\, \wt{M}_d^n \circ \bar{A}_x^n \Rightarrow 0, \,\,\mbox{ in }\,\, D(\R)\,\,\mbox{ as }\,\, n \rightarrow \infty.\]  Also, $(\mathbb{A}2)$ implies $n / \mu^n \rightarrow 1/\lambda$ as $n \rightarrow \infty$, which leads to the coefficient multiplying $W$ in \eqref{eq:WeakConvN}, and
%\begin{itemize}
%\item Standard results on renewal processes show that the delayed renewal process defined in \eqref{254} has the limiting behavior \[
%    \bar{A}_x^n \rightarrow \lambda e\,\, \mbox{ and }\,\, \wt{A}_x^n \Rightarrow \sqrt{\lambda} \sqrt{\mbox{var}(u_2)} W_1, \,\,\mbox{ in }\,\, D(\R) \,\,\mbox{ as }\,\, n \rightarrow \infty,
%    \]
%    where $W_1$ is a standard Brownian motion and $e(\cdot)$ represents the identity map, i.e., $e(t) = t$ for all $t \in \R_+$;
%    %\footnote{``$\bar{A}_x^n \rightarrow \lambda t$''?}
%\item Donsker's theorem and a random time change show
%\[
%\wt{S}^n \circ \bar{A}_x^n \Rightarrow \sqrt{\lambda} \sqrt{\mbox{var}(v_2)} W_2, \,\,\mbox{ in }\,\, D(\R) \,\,\mbox{ as }\,\, n \rightarrow \infty,
%\]
%where $W_2$ is a standard Brownian motion independent of $W_1$;
%\item Arguments similar to those in Theorem 5.1(a) in \cite{WardGlynn2005} show
%\[
%    \wt{S}_{d,x} \circ \bar{A}_x^n \Rightarrow 0 \,\,\mbox{ and }\,\, \wt{M}_d^n \circ \bar{A}_x^n \Rightarrow 0, \,\,\mbox{ in }\,\, D(\R)\,\,\mbox{ as }\,\, n \rightarrow \infty;
%\]
%\item $(\mathbb{A}2)$ implies $n / \mu^n \rightarrow 1/\lambda$ as $n \rightarrow \infty$, which leads to the coefficient multiplying $W$ in \eqref{eq:WeakConvN}.
%\end{itemize}
%Second,
$b^n \rightarrow \theta / \lambda$ as $n \rightarrow \infty$ by $(\mathbb{A}2)$. %, and modifying the arguments in Section 5.1.2 in \cite{ReedWard2008} to account for the fact that the hazard rate is not scaled shows $\wt{\epsilon}_x^n \Rightarrow 0$ as $n \rightarrow \infty$.
Finally, the linearly generalized regulator mapping in Section A.1 of \cite{WardGlynn2005} %the representation in \eqref{rep0} combined with Lemma~\ref{lem-modified}
and an application of the continuous mapping theorem establish the weak convergence in \eqref{eq:Vapprox}.
%\end{remark}
\end{proof}

Before we prove Lemma \ref{petite-set} below, we recall the definition of a petite set. Given a Markov process $\Phi=\{\Phi_t:t\geq0\}$, with the semigroup operator $P^t$, on a state space $\mathcal X$, we say that a non-empty set $A \in \mathcal B(\mathcal X)$ is $\nu_a$-petite if $\nu_a$ is a non-trivial measure on $\mathcal B(\mathcal X)$, $a(\cdot)$ is a probability measure on $(0, \infty)$ (referred to as a sampling distribution), and $K_a(x, \cdot) \geq \nu_a(\cdot)$ for all $x \in A$. See \cite{Meyn:Tweedie:1993[2]}. Here, the Markov transition function $K_a:=\int P^t a(dt)$. The set $A$ will be called simply petite when the specific measure $\nu_a$ is unimportant.

\begin{proof}[\bf Proof of Lemma \ref{petite-set}]
First, note that the set $C$ is closed by definition. To establish that it is petite, we need to find a non-trivial measure $\nu_a$ on $(\mathbb S, \mathcal B(\mathbb S))$ and a sampling distribution $a(\cdot)$ on $(0,\infty)$ such that $K_a(x, \cdot) \geq \nu_a(\cdot)$ for all $x \in C$.  We fix the scaling parameter $n\geq1$ and omit the associated superscript $n$ in this proof. %\st{The following arguments are directly transferred from the proof of Lemma 3.7 in Meyn and Down \cite{MeynDown1994}.}
{The following arguments are directly transferred from the proof of Proposition 4.8 in Bramson \cite{Bramson:2008}.} %Fix $x=(\tau, v)\in C$ and denote $\mathbf 0:=(0,0)\in \mathbb S$. Let $T_t(x,\mathbf 0):=P_x\{ \{V(t) = 0\}\cap \{A(t)=0\} \}$; that is, $T_t(x,\mathbf 0)$ is the probability the system is empty at time $t$ and no jumps have occurred in $[0,t]$.  Then, $T_t(x,\mathbf 0)=\bar F_A(t)$ when $v\leq t<\tau$ (i.e., the probability of no arrivals up until time $t$, when $t$ is greater than the initial workload $v$ but less than the initial forward recurrence time $\tau$). On the other hand, when $v\geq \tau$ and $t\geq v-\tau$, we get $T_t(x,\mathbf 0)= \bar F_A(t)\times F_S(t-(v-\tau))\bar F_D(v-\tau)$ (i.e., the initial forward recurrence time $\tau$ is expired and the arrival joins the non-empty queue with service time less than $t-(v-\tau)$). Other than these two cases, we have $T_t(x,\mathbf 0)=0$. Define $T^a(x,\{\mathbf 0\}):=\inte{0}{\infty} T_t(x,\mathbf 0)a(dt)$, $T^a(x,\mathbb S\backslash\{\mathbf 0\}):=0$, where $a$ is given as in $(\mathbb A3)$. Then, $T^a(x,\mathbb S)>0$ (i.e., nontrivial measure on $\mathbb S$) and moreover, $K_a(x,\{\mathbf 0\})\geq T^a(x,\{\mathbf 0\})$. Now, for  $A\in\mathcal B(\mathbb S)$,  $K_a(x,A)\geq \delta_0\{A\} T^a(x,\{\mathbf 0\}),$ where $\delta_0$ is the unit mass concentrated on $\mathbf 0\in \mathbb S$. Letting $\nu_a(\cdot):=\delta_0\{\cdot\}T^a(x,\{\mathbf 0\})$ completes the proof.  %For arbitrary $a\pos$, let $R_a$ denote the bounded open rectangle in $\mathbb S$ defined as $R_a\equiv [0,a)\times [0,a)$.
%\st{Fix $\kappa>0$ and consider the transition probability $P^t(x,B_1\times \{0\})$ for $B_1\in\mathcal B(\R_+)$.}
{Fix $x=(\tau, v)\in C$ and consider the transition probability $P^t(x,[0,1]\times \{0\})=\PP(\tau(t)\in [0,1], V_x(t)=0)$ for all $|x|\leq \kappa$, with $\kappa>0$ fixed.  %Choose $\bar t>0$ large enough and $\epsilon>0$ (both independent of $x$) such that there exist large $t_1(x), t_2(x)\leq  \bar t$ so that \beq\label{cl1}\PP_x(\tau(t)\in [0,1], V(t)=0)\geq \epsilon\eeq for all $t\in[t_1(x), t_2(x))$. That is, the probability is at least $\epsilon$ that the queue will be empty over some interval $[t_1(x), t_2(x))$ but not remain empty over the entire interval $[t_2(x),\bar t)$. This is owing to the fact that the interarrival time can be arbitrarily large, as assumed in $(\mathbb{A}3)$, and thus the residual arrival time $\tau(t)$ can be as small as taking values in $[0,1]$.
From $(\mathbb{A}3)$, there exists a $b\geq v+2$ so that \beqys \PP(b+1/2\leq u_{2}\leq b+1)&\geq& \epsilon_1, \\ \PP(v_{1}\leq 2)&\geq& \epsilon_2,  \eeqys for some $\epsilon_i\in (0,1), i=1,2.$  Denote the two events considered above by $G_1$ {and} $G_2$ respectively, and also define $t_1(x)\equiv \tau+b$ and $t_2(x)\equiv \tau+b+1/2$.  Suppose $\omega\in G_1\cap G_2$.
{The first arrival sees workload $(v-\tau)^+$, enters service if $d_1 > (v-\tau)^+$, and abandons the system without receiving service if $d_1 \leq (v-\tau)^+$.  Consequently, the departure time of the first arrival is
\[
t' \equiv \left\{ \begin{array}{ll} \mbox{max}(\tau,v) + v_1 &  \mbox{ if } d_1 > (v-\tau)^+,  \\ \tau+ d_1 & \mbox{ if } d_1 \leq (v-\tau)^+. \end{array} \right.
\]
Algebra shows $t' \leq \mbox{max}(\tau,v) + v_1 \leq |x| + v_1$, and so $t' \leq \tau + v + 2 \leq t_1(x)$ for $\omega \in G_2$.
}
Similarly, letting $t''$ be the arrival time of the second job, $t''=\tau+u_2\geq t_2(x)$ from the definition of $t_2(x)$ and $\omega\in G_1$. Hence, for all $t\in[t_1(x),t_2(x)]$, $V(t)=0$ and also $\tau(t)=t''-t\in[0,1]$  when $\omega\in G_1\cap G_2$. This implies that for all $t\in[t_1(x), t_2(x)]$, \beqs\PP(\tau(t)\in [0,1], V_x(t)=0) \geq  \PP(G_1\cap G_2) = \PP(G_1)\PP(G_2)\geq \epsilon_1\epsilon_2 \equiv \epsilon.\eeqs}
%(A similar result is given in Lemma 3.7 of Meyn and Down \cite{MeynDown1994} with $\tau(t)\in[0,1]$ replaced by $\tau(t)\in B_1\in\mathcal B(\R_+)$ under additional spread-out condition on the interarrival time.)
{Next, let $\bar t\equiv \kappa+b+1$ and choose the probability measure $a(\cdot)$ to be uniform over $(0,\bar t)$ so that for $B=B_1\times B_2\in \mathcal B(\mathbb S)$, \[K_a(x,B_1\times B_2):=\inte{0}{\infty} P^t(x,B_1\times B_2) a(dt)=\inte{0}{\bar t}\frac{1}{\bar t}P^t(x, B_1\times B_2)dt.\]  Set the measure $\nu(\cdot)$ on $(\mathbb S, \mathcal B(\mathbb S))$ to be uniform over the empty states with residual interarrival time in $(0,1)$, % so that its density is $\epsilon/\bar t$ with respect to this set.
that is, \[K_a(x,B_1\times B_2)\geq \epsilon 1_{\{0\in B_1\}}m((0,1)\cap B_2)\equiv \nu(B_1\times B_2),\] where $m(\cdot)$ is the Lebesgue measure on the positive real line. Since $\nu(\cdot)$ is a nontrivial measure on  $(\mathbb S, B(\mathbb S))$ (i.e., $\nu(\mathbb S)=\nu(\R_+\times \R_+)>0$), this completes the proof.} \end{proof}

\noindent{\bf \large Acknowledgements:  }{We would like to thank Junfei Huang for a discussion of the recent paper \cite{HuangGurvich18}.
Amy Ward is supported by the William S. Fishman Faculty Research Fund at the University of Chicago Booth School of Business,
and Heng-Qing Ye is supported by the HK/RGC Grant T32-102/14N
{and NSFC Grant 71520107003}. This work was performed while Chihoon Lee was visiting the Chinese University of Hong Kong, Shenzhen, during the summers of 2018 and 2019.}
%that helped us to confirm the proof of Lemma 2.}

%\bibliographystyle{plain} % outcomment this and next line in Case 1
%%\bibliography{Lee-Ward-conv-stat-revised0925} % if more than one, comma separated
%\bibliography{Lee-Ward-conv-stat-revised-4-17-2018} % if more than one,

 \newcommand{\noop}[1]{}

\end{document}